\documentclass[12pt,reqno]{amsart}
\usepackage{amsthm, amsmath, amssymb, amsbsy, amsfonts, latexsym, euscript,amscd}
\usepackage{calrsfs}
\usepackage[margin=2.2cm]{geometry}
\usepackage{tikz}
\usetikzlibrary{arrows.meta,positioning}
\usepackage{pgfplots}
\pgfplotsset{compat=1.17}
\usepackage[colorlinks=true,linkcolor=blue,citecolor=blue,urlcolor=blue]{hyperref}

\def\leq{\leqslant}

\def\geq{\geqslant}

\def\phi{\varphi}
\def\hat{\widehat}

\def\bar{\overline}

\def\kappa{\varkappa}

\newtheorem{theorem}{Theorem}[section]
\newtheorem{mtheorem}{Theorem}

\newtheorem{lemma}[theorem]{Lemma}
\newtheorem{corollary}[theorem]{Corollary}

\theoremstyle{remark}
\newtheorem{remark}[theorem]{Remark}
\newtheorem{definition}[theorem]{Definition}
\newtheorem{example}[theorem]{Example}

\numberwithin{equation}{section}

\begin{document}

\title[Idempotent manifolds and quadratic flows]{Idempotent manifolds and quadratic flows\\ of cubic forms with radial infinity-Laplacian}

\author{Daniel J. F. Fox}
\address{Departamento de Matem\'atica Aplicada a la Ingenier\'ia Industrial, Universidad Polit\'ecnica de Madrid, Ronda de Valencia 3, 28012 Madrid, Spain}
\email{daniel.fox@upm.es}

\author{Vladimir G. Tkachev}
\address{Department of Mathematics, Link\"oping University, SE-581 83 Link\"oping, Sweden}
\email{vladimir.tkatjev@liu.se}

\subjclass[2020]{Primary 17A30; Secondary 35J94, 34A34, 37J30, 58K05}
\keywords{Nonassociative algebras, idempotents, Peirce decomposition, infinity-Laplacian, cubic forms, Morse--Bott functions, Darboux points, Kovalevskaya exponents}

\begin{abstract}
We study cubic forms $u$ on a Euclidean space whose $\infty$-Laplacian is a radial multiple of $u$, $\Delta_{\infty}u=54\,\theta\,|x|^{2}u$, through the commutative algebra $\mathbb{A}$ defined by the Hessian of $u$ (a radial $\infty$-Hsiang algebra). We prove that the nonzero idempotents of $\mathbb{A}$, equivalently the positive critical points of $u$ on the unit sphere, form a compact real analytic submanifold of a sphere whose tangent space at an idempotent $e$ is the Peirce space $\mathbb{A}_{e}(\frac12)$, and that the restriction of $u$ to the unit sphere is Morse--Bott along it. We then show that the defining identity is the conservation of a quartic energy along the quadratic flow $\dot x=x^{2}$, interpret the idempotents as Darboux points of this flow and the Peirce numbers through its Kovalevskaya exponents and first integrals, and give a resonance proof that an idempotent without the Peirce eigenvalue $-\frac12$ forces a pseudocomposition identity.
\end{abstract}

\maketitle

\setcounter{tocdepth}{1}
\tableofcontents

\section{Introduction}

\subsection{Cubic forms and the \texorpdfstring{$\infty$}{infinity}-Laplacian}\label{ssec:1.1}
Let $u$ be a cubic form on a Euclidean space $\mathbb{R}^{n}$. Its gradient $\nabla u$ is quadratic and its Hessian $\operatorname{Hess}u$ is linear in $x$, so the $\infty$-Laplacian
\begin{equation}\label{eq:inflap}
\Delta_{\infty}u=\langle (\operatorname{Hess}u)\nabla u,\nabla u\rangle
\end{equation}
is a form of degree $5$. Cubic forms have long served as test objects for nonlinear elliptic equations: Hsiang's cubic minimal cones \cite{Hsiang67}, whose algebraic constructions are developed in \cite{FoxTkachevCones}, and many of the singular and nonclassical solutions of fully nonlinear elliptic equations constructed in \cite{NTV14}, which have the form $u(x)/|x|^{1+\alpha}$, are based on cubic forms. In these constructions $u$ is far from arbitrary: it has to satisfy strong rigidity conditions, and the cubic forms that work are of Jordan or Clifford origin, such as Cartan's isoparametric cubics, which are obtained by restricting the generic norm of the Jordan algebra of $3\times3$ Hermitian matrices over $\mathbb{R}$, $\mathbb{C}$, $\mathbb{H}$ or $\mathbb{O}$ to its trace-free part.

Remarkably, closely related cubic forms appear in a quite different, harmonic-analytic context. Etingof, Kazhdan and Polishchuk \cite{EtingofKazhdanPolishchuk02} showed that the cubic forms $f$ for which the Fourier transform of an elementary function with phase $f(x)/t$ can be elementary with a phase of the same type are, essentially, the generic norms of semisimple Jordan algebras of degree $3$; see Remark~\ref{rem:EKP} for the precise relation with the present paper. This suggests that such distinguished cubic forms should be studied by a uniform method. A convenient way to handle a cubic form is to pass to an algebra. Define a commutative multiplication $(x,y)\mapsto xy$ on $\mathbb{R}^{n}$, written by juxtaposition, by the relation
\begin{equation}\label{eq:mult}
\langle xy,z\rangle =\tfrac16\,d^{3}u(x,y,z),\qquad x,y,z\in\mathbb{R}^{n},
\end{equation}
and denote by $\mathbb{A}=(\mathbb{R}^{n},\,xy)$ the resulting commutative, in general nonassociative, algebra. Since the third differential $d^{3}u$ is a symmetric trilinear form and the inner product is nondegenerate, \eqref{eq:mult} determines $xy$ uniquely, and $xy=yx$. Writing $x^{2}=xx$, the Euler identity for the cubic form $u$ yields the following dictionary between calculus and algebra:
\begin{equation}\label{eq:dictionary}
\begin{aligned}
&u(x)=\langle x,x^{2}\rangle ,\\
&\nabla u(x)=3x^{2},\\
&\operatorname{Hess}u(x)=6L(x),
\end{aligned}
\end{equation}
where
\begin{equation}\label{eq:Lx}
L(x)y=xy
\end{equation}
is the multiplication operator by $x$ on $\mathbb{A}$ (left and right multiplications coincide, since $\mathbb{A}$ is commutative). By the symmetry of $d^{3}u$, the inner product is \emph{invariant}, that is,
$$
\langle xy,z\rangle =\langle x,yz\rangle \qquad  \text{for all $x,y,z$;}
$$
equivalently, every operator $L(x)$ is self-adjoint. Conversely, every commutative algebra with an invariant inner product arises in this way from its cubic form $\langle x,x^{2}\rangle $; this correspondence is the starting point of the algebraic framework developed in \cite[Chapter~6]{NTV14}. In these terms
$$
\Delta_{\infty}u=54\langle x^{2},x^{3}\rangle ,\qquad \text{ where $x^{3}=x^{2}x$.}
$$

The benefit of this dictionary is that analytic conditions on $u$, which are nonlinear systems of partial differential equations for its coefficients, become polynomial identities for the multiplication, and these can be attacked by algebraic means: idempotents, Peirce decompositions and their fusion rules, hidden Jordan and Clifford structures. In the theory of fully nonlinear uniformly elliptic equations $F(D^{2}w)=0$, viscosity solutions are known to be $C^{1,\alpha}$ \cite{CaffarelliCabre95}, and the question of whether they are always classical was open for a long time. Nadirashvili and Vl\u{a}du\c{t} \cite{NadirashviliVladut07} answered it in the negative, and the subsequent theory, summarized in \cite{NTV14}, produces singular and nonclassical solutions of the form $w=u(x)/|x|^{1+\alpha}$ with a cubic form $u$ of a very special kind, for instance a Cartan isoparametric cubic. What makes such $w$ a solution of a uniformly elliptic equation is a two-sided control of the spectrum of $D^{2}w$, and this is governed by the spectra of the operators $L(x)$ of the algebra of $u$; this mechanism is analyzed in \cite{Tkachev19}. On the geometric side, the same algebras describe cubic minimal cones: Clifford systems give an infinite family of them \cite{Tkachev10}, Hsiang's problem for cubic minimal cones leads to a rich class of algebras \cite{FoxTkachevCones}, and the cubic solutions of the eiconal equation $|\nabla u|^{2}=9|x|^{4}$ are classified by means of Jordan algebras \cite{Tkachev14}. On the algebraic side, metrized commutative algebras with an idempotent-controlled structure have become an active subject, in particular in the theory of axial algebras \cite{HallRehrenShpectorov15,GorshkovMcInroyShumbaShpectorov25}, and the eigenvalue $\frac12$ of the operators $L(e)$, which will be central below, plays a universal role in such algebras \cite{Tkachev21}.

Very recently, cubic forms of this type have also appeared in a problem of approximation theory. In \cite{Tkachev26}, Hsiang cubics are studied as spherical analogues of Chebyshev polynomials. Such a cubic equioscillates on the unit sphere, and its maximum set consists of the normalized idempotents. In the problem of least deviation from zero on the sphere among harmonic cubics, Cartan's isoparametric cubics turn out to be the minimizers in their dimensions, and the Kolmogorov criterion for stationarity becomes the existence of a positive measure on the idempotent set, which in turn produces quadrature formulas on the sphere. In that context the geometry of the set of maximum points of $u$ on the sphere, its dimension, its tangent spaces and the behavior of $u$ near it, are exactly the information one needs. This is one of the motivations for the results of Subsection~\ref{ssec:1.3} below, where this set is described for a broad class of cubic forms.

In this paper we study cubic forms whose $\infty$-Laplacian is a radial multiple of the form itself,
\begin{equation}\label{eq:radialinf}
\Delta_{\infty}u=54\,\theta\,|x|^{2}\,u\qquad\text{for some constant }\theta\neq0 ,
\end{equation}
together with the associated algebra (the case $\theta=0$ is trivial, see Remark~\ref{rem:thetazero}). The equation \eqref{eq:radialinf} is the cubic analogue of Hsiang's equation for minimal cones \cite{Hsiang67}, with the ordinary Laplacian of $u$ removed. Its algebraic form makes sense over more general fields, but throughout this paper the base field is $\mathbb{R}$ and the metric is Euclidean (see Remark~\ref{rem:standing}). Let $\mathbb{A}$ be a finite-dimensional commutative, not necessarily associative, real algebra. A \emph{metric} on $\mathbb{A}$ is a nondegenerate symmetric bilinear form $\langle \cdot,\cdot\rangle $ which is invariant, that is, every $L(x)$ is self-adjoint; it is \emph{Euclidean} if it is positive definite. We put
\begin{equation}\label{eq:qu}
\begin{aligned}
q(x)&=\langle x,x\rangle ,\\
u(x)&=\langle x,x^{2}\rangle .
\end{aligned}
\end{equation}

\begin{definition}\label{def:hsiang}
A commutative algebra $\mathbb{A}$ with a metric is \emph{radial $\infty$-Hsiang} with constant $\theta\in\mathbb{R}^{\times}$ if
\begin{equation}\label{eq:hsiang}
\langle x^{2},x^{3}\rangle =\theta\,\langle x,x\rangle \,\langle x,x^{2}\rangle \qquad\text{for all }x\in\mathbb{A} .
\end{equation}
\end{definition}

The first observation, which motivates the definition, is that \eqref{eq:hsiang} is just the algebraic form of \eqref{eq:radialinf}. Indeed, for a Euclidean space with the multiplication \eqref{eq:mult}, the dictionary \eqref{eq:dictionary} gives $\Delta_{\infty}u=\langle 6L(x)\,3x^{2},3x^{2}\rangle =54\langle x^{2},x^{3}\rangle$ and $|x|^{2}u=\langle x,x\rangle \langle x,x^{2}\rangle$. Hence a cubic form $u$ satisfies \eqref{eq:radialinf} if and only if its algebra is radial $\infty$-Hsiang with the same constant $\theta$, and, conversely, every radial $\infty$-Hsiang algebra over $\mathbb{R}$ with a positive definite metric arises in this way from its cubic form $u(x)=\langle x,x^{2}\rangle$. Two classical families of nonassociative algebras satisfy stronger, vector-valued identities. A \emph{pseudocomposition algebra} \cite{MeybergOsborn93} satisfies
$$x^{3}=b(x,x)x
$$ for \emph{some} symmetric bilinear form $b$, and an \emph{admissible cubic algebra} \cite{ElduqueOkubo00} satisfies
$$
x^{2}x^{2}=N(x)x
$$ for \emph{some} cubic form $N$ (the letter $N$ is used only here). The word \emph{radial} means that these forms are tied to the given metric:
\begin{equation}\label{eq:radialclasses}
\begin{aligned}
x^{3}&=\theta\langle x,x\rangle \,x &&\text{(radial pseudocomposition)},\\
x^{2}x^{2}&=\theta\langle x,x^{2}\rangle \,x &&\text{(radial admissible cubic)}.
\end{aligned}
\end{equation}
The word \emph{radial} has the same meaning in Definition~\ref{def:hsiang}, where both sides are expressed through the one metric. Put
\begin{equation}\label{eq:defects}
\begin{aligned}
\mathcal{R}(x)&=x^{3}-\theta q(x)\,x,\\
\mathcal{S}(x)&=x^{2}x^{2}-\theta u(x)\,x .
\end{aligned}
\end{equation}
Since $\langle x^{2}x^{2},x\rangle =\langle x^{2},x^{3}\rangle $ by invariance,
\begin{equation}\label{eq:shadow}
\langle x^{2},x^{3}\rangle -\theta q(x)u(x)=\langle \mathcal{R}(x),x^{2}\rangle =\langle \mathcal{S}(x),x\rangle .
\end{equation}
Thus \eqref{eq:hsiang} is a common scalar consequence of the two vector identities $\mathcal{R}=0$ and $\mathcal{S}=0$, and both families in \eqref{eq:radialclasses} are radial $\infty$-Hsiang. We call $\mathcal{R}$ and $\mathcal{S}$ the \emph{pseudocomposition defect} and the \emph{admissible defect}.

\subsection{Idempotents and Peirce numbers}\label{ssec:1.2}
An element $e\in\mathbb{A}$ with $e^{2}=e$ is an \emph{idempotent}, and an element $x$ with $x^{2}=0$ is a \emph{$2$-nilpotent}. We denote by $\operatorname{Idm}(\mathbb{A})$ the set of nonzero idempotents. Analytically, idempotents have two familiar meanings. First, the rescaled idempotents are exactly the positive critical points of $u$ on the unit sphere (Lemma~\ref{lem:critical}), since the Lagrange equation $\nabla u=\lambda x$ reads $x^{2}=\frac{\lambda}{3}x$. And for the quadratic system $\dot x=x^{2}$, every idempotent produces the explicit solution $x(t)=-t^{-1}e$, which runs along a straight line (Subsection~\ref{ssec:1.4}).

In a radial $\infty$-Hsiang algebra every $e\in\operatorname{Idm}(\mathbb{A})$ satisfies $\theta\langle e,e\rangle =1$, and the self-adjoint operator $L(e)$ is diagonalizable with eigenvalue $1$ on $\mathbb{R} e$ and all other eigenvalues in $\{-1,-\frac12,\frac12\}$ (Lemma~\ref{prop:peirce}). We write
\begin{equation}\label{eq:peirceno}
\mathbb{A}_{e}(\mu)=\ker(L(e)-\mu\,\operatorname{id}),\qquad n_{\mu}(e)=\dim\mathbb{A}_{e}(\mu),\qquad \mu\in\{-1,-\tfrac12,\tfrac12\},
\end{equation}
and call the $n_{\mu}(e)$ the \emph{Peirce numbers} of $e$.

\begin{definition}\label{def:complexity}
The number of distinct eigenvalues of $L(e)$, that is, $1$ plus the number of nonzero Peirce numbers of $e$, is called the \emph{Peirce complexity} of $e$. It is at most $4$, and we say that $e$ has \emph{maximal} (respectively, \emph{submaximal}) Peirce complexity if it equals $4$ (respectively, $3$).
\end{definition}

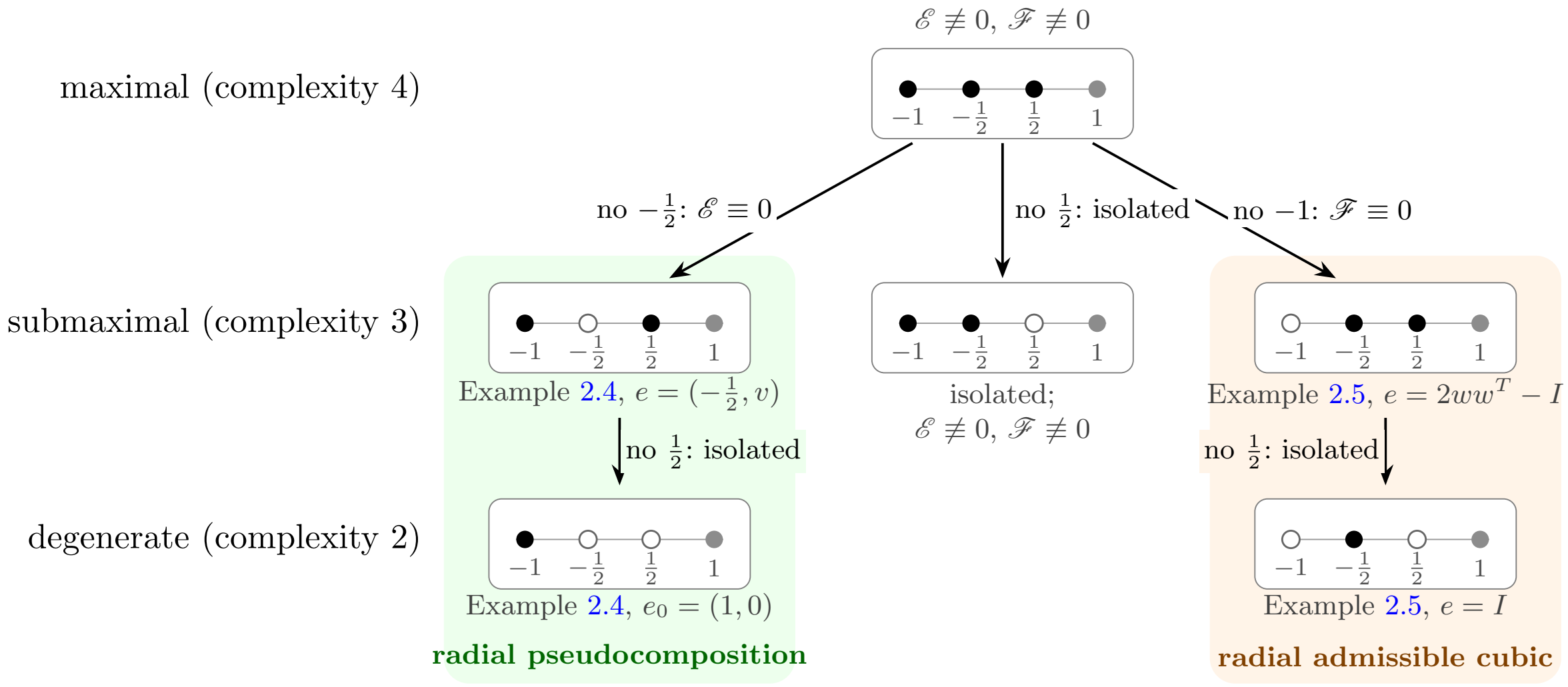
\begin{figure}[ht]
\centering
\resizebox{0.92\textwidth}{!}{%
\begin{tikzpicture}[>=Stealth,
  on/.style={circle,fill=black,inner sep=0pt,minimum size=5.5pt},
  off/.style={circle,draw=black!60,fill=white,inner sep=0pt,minimum size=5.5pt,line width=0.6pt},
    pics/spec/.style n args={3}{code={
    \fill[rounded corners=4pt,white] (-1.45,-0.55) rectangle (1.45,0.45);
    \draw[rounded corners=4pt,black!50] (-1.45,-0.55) rectangle (1.45,0.45);
    \draw[black!40] (-1.15,0)--(1.15,0);
    \node[#1] at (-1.05,0) {}; \node[#2] at (-0.35,0) {}; \node[#3] at (0.35,0) {}; \node[on,fill=black!45] at (1.05,0) {};
    \node[font=\tiny,text=black!70] at (-1.05,-0.32) {$-1$};
    \node[font=\tiny,text=black!70] at (-0.35,-0.32) {$-\frac12$};
    \node[font=\tiny,text=black!70] at (0.35,-0.32) {$\frac12$};
    \node[font=\tiny,text=black!70] at (1.05,-0.32) {$1$};
  }},
  lbl/.style={font=\scriptsize,align=center,fill=white,inner sep=1.5pt},
  ex/.style={font=\scriptsize,align=center,text=black!75}]
\fill[green!7,rounded corners=8pt] (-6.2,-1.85) rectangle (-2.3,-6.6);
\fill[orange!9,rounded corners=8pt] (2.3,-1.85) rectangle (6.2,-6.6);
\node[font=\scriptsize\bfseries,text=green!40!black] at (-4.25,-6.3) {radial pseudocomposition};
\node[font=\scriptsize\bfseries,text=orange!50!black] at (4.25,-6.3) {radial admissible cubic};
\node[font=\small,anchor=east] at (-6.3,0) {maximal (complexity $4$)};
\node[font=\small,anchor=east] at (-6.3,-2.6) {submaximal (complexity $3$)};
\node[font=\small,anchor=east] at (-6.3,-5.0) {degenerate (complexity $2$)};
\pic at (0,0) {spec={on}{on}{on}};
\pic at (-4.25,-2.6) {spec={on}{off}{on}};
\pic at (4.25,-2.6) {spec={off}{on}{on}};
\pic at (-4.25,-5.0) {spec={on}{off}{off}};
\pic at (4.25,-5.0) {spec={off}{on}{off}};
\node[ex] at (0,0.75) {$\mathcal{E}\not\equiv0$, $\mathcal{F}\not\equiv0$};
\draw[->,thick] (-1.0,-0.6) -- node[lbl,left=4pt] {no $-\frac12$:\ $\mathcal{E}\equiv0$} (-3.7,-2.1);
\draw[->,thick] (1.0,-0.6) -- node[lbl,right=4pt] {no $-1$:\ $\mathcal{F}\equiv0$} (3.7,-2.1);
\draw[->,thick] (-4.25,-3.65) -- node[lbl,right,fill=green!7] {no $\frac12$:\ isolated} (-4.25,-4.4);
\draw[->,thick] (4.25,-3.65) -- node[lbl,left,fill=orange!9] {no $\frac12$:\ isolated} (4.25,-4.4);
\node[ex] at (-4.25,-5.75) {Example~\ref{ex:rpc}, $e_{0}=(1,0)$};
\node[ex] at (4.25,-5.75) {Example~\ref{ex:rac}, $e=I$};
\node[ex] at (-4.25,-3.38) {Example~\ref{ex:rpc}, $e=(-\frac12,v)$};
\node[ex] at (4.25,-3.38) {Example~\ref{ex:rac}, $e=2ww^{T}-I$};
\pic at (0,-2.6) {spec={on}{on}{off}};
\draw[->,thick] (0,-0.6) -- node[lbl,right=2pt] {no $\frac12$:\ isolated} (0,-2.1);
\node[ex] at (0,-3.6) {isolated;\\ $\mathcal{E}\not\equiv0$, $\mathcal{F}\not\equiv0$};
\end{tikzpicture}
}
\caption{Peirce complexity of an idempotent $e$ of a radial $\infty$-Hsiang algebra. Each box shows the spectrum of $L(e)$: a filled dot means that the eigenvalue occurs, an open dot that it does not; the eigenvalue $1$ always occurs. Each arrow removes one eigenvalue on $e^{\perp}$ and shows the consequence: without $-\frac12$ the energy vanishes, $\mathcal{E}\equiv0$ (Theorem~\ref{mthm:C}); without $-1$ the second integral vanishes, $\mathcal{F}\equiv0$ (the rigidity of the eigenvalue $-1$; Corollary~\ref{cor:rac2} in the degenerate case); without $\frac12$ the idempotent is isolated (Theorem~\ref{mthm:A}).}
\label{fig:complexity}
\end{figure}

 The possible patterns, and what they mean for the algebra and for the objects studied below, are shown in Figure~\ref{fig:complexity}.

 Linearizing the defects at $e$ (Lemma~\ref{lem:defects}) shows that pseudocomposition algebras allow only the eigenvalues $\frac12$ and $-1$ on $e^{\perp}$, and admissible cubic algebras only $\frac12$ and $-\frac12$; radial $\infty$-Hsiang algebras allow the union.

In a forthcoming paper, we show by purely algebraic means (fusion rules for the Peirce spaces and a Jordan algebra structure on $\mathbb{R} e\oplus\mathbb{A}_{e}(-\frac12)$), which apply over any field of characteristic not dividing $30$, that the absence of one of the two eigenvalues at a \emph{single} idempotent propagates to the whole algebra: if $n_{-1/2}(e)=0$ for some $e\in\operatorname{Idm}(\mathbb{A})$, then $\mathbb{A}$ is radial pseudocomposition, and if $n_{-1}(e)=0$ for some $e\in\operatorname{Idm}(\mathbb{A})$, then $\mathbb{A}$ is radial admissible cubic. In particular, whether the eigenvalue $-1$ (respectively, $-\frac12$) occurs does not depend on the idempotent, while the occurrence of the eigenvalue $\frac12$ and all multiplicities may depend on it (Examples~\ref{ex:rpc} and \ref{ex:rac}). None of the results of this forthcoming paper is used in the proofs below. We shall refer to these two statements as the \emph{rigidity of the eigenvalue $-\frac12$} and the \emph{rigidity of the eigenvalue $-1$}.

The present paper has two aims, both of an analytic and geometric nature. The first is to describe $\operatorname{Idm}(\mathbb{A})$ as a geometric object. The second is to understand the two rigidity statements through the dynamics of $\dot x=x^{2}$, in particular through its first integrals and its Kovalevskaya exponents. In this way we obtain a self-contained proof of the rigidity of the eigenvalue $-\frac12$, which explains why the absence of this eigenvalue at one idempotent controls the whole algebra, and we identify precisely what obstructs the same mechanism for the eigenvalue $-1$.

\subsection{The idempotent manifold}\label{ssec:1.3}
Our first result describes the set of all idempotents, which acquires a transparent geometric meaning. To appreciate it, recall what one should expect in general. The nonzero idempotents of an $n$-dimensional commutative algebra are the nonzero solutions of the system of $n$ quadratic equations $x^{2}=x$. For a generic algebra this system has only isolated solutions: over $\mathbb{C}$, by B\'ezout's theorem, a generic algebra has exactly $2^{n}-1$ nonzero idempotents, and a generic real algebra has finitely many of them \cite{KrasnovTkachev19}. An idempotent $e$ can fail to be isolated only if the linearization $2L(e)-\operatorname{id}$ of the equation $x^{2}=x$ is degenerate, that is, only if $\frac12$ is an eigenvalue of $L(e)$, which is a nongeneric condition. When this happens, nothing in general prevents the idempotents near $e$ from forming a singular set, or from consisting of $e$ alone: the solutions $y$ of the linearized equation $2ey=y$ need not be tangent to any curve of idempotents. In a radial $\infty$-Hsiang algebra, on the other hand, the eigenvalue $\frac12$ is the rule rather than the exception, and the following theorem shows that the idempotents nevertheless behave in the best possible way. This is a very special feature of these algebras.

Under the standing assumption, $\mathbb{A}$ is a real algebra with a Euclidean metric, and we call such an algebra \emph{Euclidean}. We keep this adjective in the statements whose proofs use the positivity of the metric in an essential way.

\begin{mtheorem}\label{mthm:A}
Let $\mathbb{A}$ be a Euclidean radial $\infty$-Hsiang algebra with nonzero multiplication. Then $\theta>0$, the set $\operatorname{Idm}(\mathbb{A})$ is nonempty, and it is a compact real analytic submanifold (whose connected components may have different dimensions) of the sphere
\begin{equation}\label{eq:sphere}
\{x\in\mathbb{A}:\ \theta\langle x,x\rangle =1\}.
\end{equation}
Its tangent space at $e\in\operatorname{Idm}(\mathbb{A})$ is $\mathbb{A}_{e}(\frac12)$. In particular the connected component of $\operatorname{Idm}(\mathbb{A})$ through $e$ has dimension $n_{1/2}(e)$, and all three Peirce numbers are constant on each connected component.
\end{mtheorem}

Different components can have different dimensions: in Example~\ref{ex:rpc}, $\operatorname{Idm}(\mathbb{A})$ is the union of a point and a round sphere, and in Example~\ref{ex:rac} it is the union of a point and a real projective plane. For exact Hsiang algebras (those with $\operatorname{tr}L(x)\equiv0$ and satisfying Hsiang's identity of \cite{NTV14}), the analogous statement, with a single type of component, is proved in a monograph of the second author on Hsiang algebras, which is in preparation; it is not used here. The proof given here follows the same strategy, a Lyapunov--Schmidt reduction of the equation $x^{2}=x$ (which in our finite-dimensional setting amounts to an application of the implicit function theorem), but it rests on an identity that we derive for radial $\infty$-Hsiang algebras (Lemma~\ref{lem:defectidentity}). This identity shows that the bifurcation equations of the reduction hold automatically.

The idempotent manifold also has a variational meaning. Let $\Sigma=\{x:\langle x,x\rangle =1\}$ be the unit sphere.

\begin{mtheorem}\label{mthm:MB}
Let $\mathbb{A}$ be a Euclidean radial $\infty$-Hsiang algebra with nonzero multiplication.
\begin{itemize}
\item[{(i)}] The critical points of $u|_{\Sigma}$ with positive critical value are exactly the points of $\sqrt{\theta}\,\operatorname{Idm}(\mathbb{A})$. They all have the critical value $\sqrt{\theta}$, which is the maximum of $u$ on $\Sigma$; so every positive critical point is a global maximum.
\item[{(ii)}] $u|_{\Sigma}$ is a Morse--Bott function. Its critical set is the disjoint union of three compact real analytic submanifolds:
\begin{itemize}
\item
the maximum set $\sqrt{\theta}\,\operatorname{Idm}(\mathbb{A})$,
\item the minimum set $-\sqrt{\theta}\,\operatorname{Idm}(\mathbb{A})$, and
\item the set $\operatorname{Nil}_{\Sigma}(\mathbb{A})=\{x\in\Sigma:x^{2}=0\}$ of normalized $2$-nilpotents, on which $u=0$.
\end{itemize}
At $\sqrt{\theta}e$ the Hessian of $u|_{\Sigma}$ is $-3\sqrt{\theta}\,\mathcal{K}_{e}=3\sqrt{\theta}\,(2L(e)-\operatorname{id})$ on $e^{\perp}$, which vanishes exactly on the tangent space $\mathbb{A}_{e}(\frac12)$ and is negative definite on $\mathbb{A}_{e}(-1)\oplus\mathbb{A}_{e}(-\frac12)$; at $x_{0}\in\operatorname{Nil}_{\Sigma}(\mathbb{A})$ it is $6L(x_{0})$ on $x_{0}^{\perp}$, whose kernel $\ker L(x_{0})\cap x_{0}^{\perp}$ is the tangent space of $\operatorname{Nil}_{\Sigma}(\mathbb{A})$.
\item[{(iii)}] The connected components of $\{x\in \Sigma:u(x)>0\}$ are in bijection with the connected components of $\operatorname{Idm}(\mathbb{A})$, each component of $\{u>0\}$ containing exactly one component of $\sqrt{\theta}\,\operatorname{Idm}(\mathbb{A})$.
\end{itemize}
\end{mtheorem}

\subsection{The quadratic flow, Darboux points and Kovalevskaya exponents}\label{ssec:1.4}
The second group of results is algebraic in nature. It concerns the quadratic differential equation
\begin{equation}\label{eq:flow}
\dot x=x^{2}
\end{equation}
attached to $\mathbb{A}$. The correspondence between quadratic systems and commutative algebras goes back to Markus \cite{Markus60} and was developed by R\"ohrl \cite{Rohrl77}, Walcher \cite{Walcher91} and Kinyon and Sagle \cite{KinyonSagle95}; see \cite{KrasnovTkachev19,Krasnov23} for its relation to idempotents and for applications to differential equations. By \eqref{eq:dictionary}, \eqref{eq:flow} is the gradient flow of $\frac13u$. The same quadratic map $x\mapsto x^{2}=\frac13\nabla u(x)$ also generates the multicomponent system of hydrodynamic type $x_{t}=(x^{2})_{s}=2L(x)\,x_{s}$, whose characteristic velocities are the eigenvalues of $2L(x)$; for the Cartan--M\"unzner polynomials of isoparametric hypersurfaces, Ferapontov \cite{Ferapontov95} showed that such systems, restricted to homogeneous isoparametric hypersurfaces, are nondiagonalizable but integrable. We use two quantities built from the metric,
\begin{equation}\label{eq:EFintro}
\begin{aligned}
\mathcal{E}(x)&=\langle x^{2},x^{2}\rangle -\theta\langle x,x\rangle ^{2},\\
\mathcal{F}(x)&=\langle x^{2},x^{2}x^{2}\rangle -\theta\langle x,x^{2}\rangle ^{2}.
\end{aligned}
\end{equation}
Along a solution of \eqref{eq:flow}, $\mathcal{E}=|\dot x|^{2}-\theta|x|^{4}$, and we call $\mathcal{E}$ the \emph{energy}.

\begin{mtheorem}\label{mthm:B}
Let $\mathbb{A}$ be a commutative nonassociative algebra with a metric and let $\theta\in\mathbb{R}^{\times}$.
\begin{itemize}
\item[{(i)}] $\mathbb{A}$ is radial $\infty$-Hsiang with constant $\theta$ if and only if $\mathcal{E}$ is a first integral of \eqref{eq:flow}; in that case $\mathcal{F}$ is a first integral as well.
\item[{(ii)}] $\nabla \mathcal{E}=4\mathcal{R}$ and $\nabla \mathcal{F}=6L(x)\mathcal{S}(x)$, where $\mathcal{R}$ and $\mathcal{S}$ are the defects \eqref{eq:defects}. Hence $\mathcal{E}\equiv0$ if and only if $\mathbb{A}$ is radial pseudocomposition, and, if $\mathbb{A}$ is radial $\infty$-Hsiang with a nonzero idempotent, $\mathcal{F}\equiv0$ if and only if $\mathbb{A}$ is radial admissible cubic.
\item[{(iii)}] If $\mathbb{A}$ is radial $\infty$-Hsiang and $e\in\operatorname{Idm}(\mathbb{A})$, then $\mathcal{E}(e)=\mathcal{F}(e)=0$ and $\nabla \mathcal{E}(e)=\nabla \mathcal{F}(e)=0$, and the restrictions of their Hessians at $e$ to $e^{\perp}$ have ranks $n_{-1/2}(e)$ and $n_{-1}(e)$, respectively.
\end{itemize}
\end{mtheorem}

Part (iii) needs a word of explanation. Since $\mathcal{E}$ and $\mathcal{F}$ vanish at every idempotent $e$ together with their gradients, the leading term of each of them near $e$ is its Hessian, a quadratic form on $e^{\perp}$. We shall see (Theorem~\ref{thm:hessians}) that the Hessian of $\mathcal{E}$ at $e$ is nondegenerate on $\mathbb{A}_{e}(-\frac12)$ and vanishes on the other Peirce spaces, while the Hessian of $\mathcal{F}$ is nondegenerate on $\mathbb{A}_{e}(-1)$ and vanishes on the others. Thus $n_{-1/2}(e)=0$ means exactly that the Hessian of $\mathcal{E}$ at $e$ vanishes, and $n_{-1}(e)=0$ that the Hessian of $\mathcal{F}$ at $e$ vanishes. Combined with (ii), the two rigidity statements of Subsection~\ref{ssec:1.2} can therefore be read as follows: if the Hessian of one of the two global conserved quantities $\mathcal{E}$, $\mathcal{F}$ vanishes at a single idempotent, then this quantity vanishes identically. To explain why this should be so, we use a classical tool from the theory of integrability of homogeneous differential equations.

Let $c\neq0$ satisfy $c^{2}=-c$. Then $x(t)=t^{-1}c$ solves \eqref{eq:flow}, because both sides equal $-t^{-2}c$; this solution moves along the straight line $\mathbb{R} c$. Such a point $c$ (or the line through it) is called a \emph{Darboux point}. For homogeneous systems, Darboux points play the role that equilibria play in the local theory: they are the simplest explicit solutions, and one studies the system by linearizing along them. For \eqref{eq:flow} the variational equation along $x(t)=t^{-1}c$ is
\begin{equation}\label{eq:variational}
\dot\xi=2x(t)\,\xi=\frac{2}{t}\,L(c)\,\xi .
\end{equation}
This is an Euler equation, and its solutions are combinations of the powers $t^{\rho-1}\xi_{0}$, where $\rho$ runs over the eigenvalues of the \emph{Kovalevskaya matrix}
\begin{equation}\label{eq:Kc}
\mathcal{K}(c)=\operatorname{id}+2L(c)
\end{equation}
and $\xi_{0}$ over the corresponding eigenvectors. These eigenvalues are the \emph{Kovalevskaya exponents} of $c$ \cite{Yoshida83}. The name recalls S.~Kovalevskaya's discovery of her integrable top \cite{Kowalevski89}. She looked for the cases in which the general solution is meromorphic in complex time, and for this the exponents must, in particular, be integers. This idea has been developed into a systematic theory of necessary conditions for integrability by Yoshida \cite{Yoshida83}, Kozlov \cite{Kozlov96}, Furta \cite{Furta96} and, through differential Galois theory, by Morales-Ruiz and Ramis \cite{MoralesRuiz99}; see also the book \cite{Goriely01}. In the Hamiltonian setting, Darboux points of homogeneous potentials and the eigenvalues of the Hessian of the potential at them have been studied in great detail by Maciejewski and Przybylska \cite{MaciejewskiPrzybylska04,MaciejewskiPrzybylska05,Przybylska09}.

For us the relevant fact is the most elementary one, due to Yoshida \cite{Yoshida83} and Kozlov \cite{Kozlov96}. If a homogeneous first integral of degree $d$ is expanded near a Darboux point, then its lowest-order part can contain only \emph{resonant} monomials, namely those whose exponents add up to $d$ (Lemma~\ref{lem:resonance}). Every nonzero idempotent gives the Darboux point $c=-e$, and then, by \eqref{eq:Kc}, $\mathcal{K}(-e)=\operatorname{id}-2L(e)$, an operator that we denote by $\mathcal{K}_{e}$ (see \eqref{eq:Ke}). So the Peirce decomposition is precisely the eigenspace decomposition of the Kovalevskaya matrix:
\begin{equation}\label{eq:kovalevskaya}
\rho=1-2\mu\in\{-1,\ 3,\ 2,\ 0\}\qquad\text{for}\qquad \mu\in\{1,\ -1,\ -\tfrac12,\ \tfrac12\}.
\end{equation}
The exponent $\rho=-1$ on $\mathbb{R} e$ is always present and corresponds to the time shift $t\mapsto t-t_{0}$ of the special solution. The exponent $\rho=0$ gives solutions $\xi=t^{-1}\xi_{0}$ of \eqref{eq:variational} that decay at the same rate as $x(t)$ itself, that is, infinitesimal motions to neighboring straight-line solutions. By Theorem~\ref{mthm:A}, these are exactly the directions tangent to the idempotent manifold. The exponents $3$ and $2$ correspond to the transversal directions, and Theorem~\ref{mthm:B}(iii) shows that they are detected by the Hessians of $\mathcal{F}$ and $\mathcal{E}$.

Let us now count resonances. We write $\rho_{\mu}=1-2\mu$ for the exponent of the Peirce space $\mathbb{A}_{e}(\mu)$, so that $\rho_{-1}=3$, $\rho_{-1/2}=2$ and $\rho_{1/2}=0$.

\begin{definition}
Choose a basis $b_{1},\dots,b_{N}$ of $e^{\perp}$ consisting of eigenvectors of $L(e)$, say $L(e)b_{i}=\mu_{i}b_{i}$, and write $z=\sum_{i}z_{i}b_{i}$ for $z\in e^{\perp}$; we call $z_{1},\dots,z_{N}$ \emph{Peirce coordinates} at $e$. The coordinate $z_{i}$ is given the \emph{weight} $\rho_{i}=1-2\mu_{i}$, and the weight of a monomial $z^{\alpha}=z_{1}^{\alpha_{1}}\cdots z_{N}^{\alpha_{N}}$ is $\alpha_{1}\rho_{1}+\dots+\alpha_{N}\rho_{N}$, the sum of the weights of its factors. A monomial is \emph{resonant} for a homogeneous first integral of degree $d$ if its weight equals $d$.
\end{definition}
We record a resonance by listing the nonzero weights of its factors. For instance, $2+2+2=6$ stands for a monomial with three factors from $\mathbb{A}_{e}(-\frac12)$, that is, $\rho_{-1/2}+\rho_{-1/2}+\rho_{-1/2}=6$, together with any number of factors from $\mathbb{A}_{e}(\frac12)$, which have weight $0$; such a monomial is resonant for an integral of degree $6$. For the energy $\mathcal{E}$ in \eqref{eq:EFintro}, of degree $4$, the only resonance is $2+2=4$, since $4$ is not a sum of copies of $3$. If $\mathbb{A}_{e}(-\frac12)=0$, there is no resonant monomial at all, and $\mathcal{E}$ must vanish. This gives the following result.

\begin{mtheorem}\label{mthm:C}
Let $\mathbb{A}$ be radial $\infty$-Hsiang. If some $e\in\operatorname{Idm}(\mathbb{A})$ has $n_{-1/2}(e)=0$, then the energy $\mathcal{E}$ in \eqref{eq:EFintro} vanishes identically. Consequently $\mathbb{A}$ is radial pseudocomposition and $n_{-1/2}(e')=0$ for all $e'\in\operatorname{Idm}(\mathbb{A})$.
\end{mtheorem}

\begin{figure}[ht]
\centering
\resizebox{0.85\textwidth}{!}{%
\begin{tikzpicture}[scale=0.8,transform shape,
  box/.style={draw,rounded corners=3pt,minimum width=3.0cm,text width=2.9cm,minimum height=0.8cm,align=center,font=\small},
  >=Stealth]
\foreach \x/\m/\r/\t [count=\i] in {
   0/1/{-1}/{the line $\mathbb{R} e$\\(time shift)},
   3.9/{\tfrac12}/0/{tangent to\\$\operatorname{Idm}(\mathbb{A})$},
   7.8/{-\tfrac12}/2/{seen by $\mathcal{E}$\\(resonance $2+2=4$)},
   11.7/{-1}/3/{seen by $\mathcal{F}$\\(resonance $3+3=6$)}}{
  \node[box,fill=gray!8] (m\i) at (\x,2.4) {$\mu=\m$};
  \node[box,fill=blue!6] (r\i) at (\x,1.2) {$\rho=\r$};
  \node[box,minimum height=1.1cm] (t\i) at (\x,-0.1) {\t};
  \draw[->] (m\i) -- (r\i);
  \draw[->] (r\i) -- (t\i);
}
\end{tikzpicture}
}
\caption{Peirce eigenvalues $\mu$ of an idempotent $e$, the Kovalevskaya exponents $\rho$ of $\dot x=x^{2}$ at the Darboux point $-e$, and their geometric and dynamical meaning.}
\label{fig:main}
\end{figure}
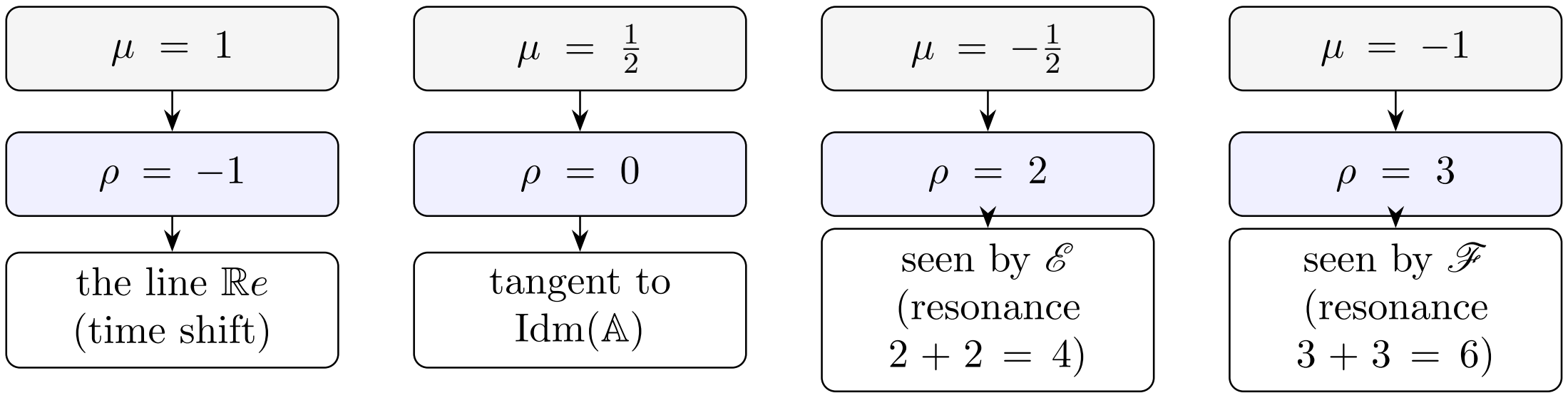

This provides an alternative (dynamical) proof of the rigidity of the eigenvalue $-\frac12$. For the second integral $\mathcal{F}$ in \eqref{eq:EFintro}, of degree $6$, besides the resonance $3+3=6$ there is also the resonance $2+2+2=6$, and this is exactly why the rigidity of the eigenvalue $-1$ does not follow in the same way. We identify the obstruction explicitly (Lemma~\ref{prop:racobstruction}) and show that it disappears when, in addition, $n_{1/2}(e)=0$ (Corollary~\ref{cor:rac2}). Figure~\ref{fig:main} summarizes the dictionary between the Peirce data and the dynamics.

\begin{remark}
Kovalevskaya exponents thus serve in this paper mainly as a language: they explain which Peirce data a first integral can detect, and the only technical input from the theory of integrability is the resonance lemma, which we prove in full. Deeper results enter only in one place. Combining Theorem~\ref{mthm:A} with the Morales-Ruiz--Ramis theorem in the form given in \cite{MaciejewskiPrzybylska05}, we observe that the Newtonian system $\ddot z=z^{2}$ attached to a Euclidean radial $\infty$-Hsiang algebra of dimension at least $2$ is never integrable with meromorphic first integrals (Lemma~\ref{prop:newton}). This second-order system should not be confused with the first-order system \eqref{eq:flow}, whose Kovalevskaya exponents \eqref{eq:kovalevskaya} are all integers, so that Kovalevskaya's test gives no obstruction for it.
\end{remark}

\subsection{Relation to the algebraic approach}\label{ssec:1.5}
In the forthcoming paper mentioned in Subsection~\ref{ssec:1.2}, the structure theory of radial $\infty$-Hsiang algebras is developed by nonassociative methods, and the algebras with an idempotent of submaximal Peirce complexity are classified. No result of that paper is used here. The present paper is complementary and is addressed to readers interested in geometry, analysis and dynamics. Theorems~\ref{mthm:A} and \ref{mthm:MB} and the description of Peirce numbers through first integrals are new. Some algebraic facts are proved again, by different means and independently of that paper: the Peirce decomposition (Lemma~\ref{prop:peirce}), which we derive from a single gradient identity; the rigidity of the eigenvalue $-\frac12$, which here follows from the conservation of energy and a resonance count (Theorem~\ref{mthm:C}); and the special case $n_{-1}(e)=n_{1/2}(e)=0$ of the rigidity of the eigenvalue $-1$ (Corollary~\ref{cor:rac2}), for which the algebraic approach gives in addition the isomorphism type of the algebra. We do not reprove the rigidity of the eigenvalue $-1$ in general; Section~\ref{sec:rac} explains, in terms of resonances, what makes it harder.

\begin{remark}[Standing assumption]\label{rem:standing}
Throughout the paper the base field is $\mathbb{R}$ and the metric is Euclidean, that is, positive definite. Several of the arguments below are formal. The gradient identity of Lemma~\ref{lem:gradient}, the Peirce decomposition, the first integrals of Section~\ref{sec:flow} and the resonance lemma are polynomial identities, and their proofs remain valid over more general fields, apart from some small characteristics. The restriction is nevertheless the natural setting for the present paper. Its subject is the correspondence between algebraic notions, such as idempotents, Peirce numbers and the identities $\mathcal{R}=0$ and $\mathcal{S}=0$, and analytic and dynamical notions, such as critical points of $u$ on the sphere, gradient flows, Darboux points, Kovalevskaya exponents, resonances and first integrals. The dynamical notions have their direct meaning over $\mathbb{R}$. The derivation $\partial$ of Section~\ref{sec:flow} is the derivative along actual solutions of \eqref{eq:flow}, a first integral is a conserved quantity, a Darboux point yields the explicit solution $x(t)=-t^{-1}e$, and the Kovalevskaya exponents are the exponents of the solutions $t^{\rho-1}\xi_{0}$ of the variational equation \eqref{eq:variational}.

The Euclidean metric plays a further role. It identifies \eqref{eq:flow} with the gradient flow of $\frac13u$, and hence identifies \eqref{eq:hsiang} with the conservation of the energy $|\dot x|^{2}-\theta|x|^{4}$, which is the link with the $\infty$-Laplacian. Positivity also guarantees the existence of idempotents and the compactness of $\operatorname{Idm}(\mathbb{A})$ (Lemma~\ref{prop:existence}), underlies the Morse--Bott description of Theorem~\ref{mthm:MB}, and makes the case $\theta=0$ trivial (Remark~\ref{rem:thetazero}).

Finally, the resonance mechanism itself is characteristic-sensitive. Over a field of characteristic $p>0$, the argument of Lemma~\ref{lem:resonance} shows only that every monomial $z^{\alpha}$ in the lowest-order part of a homogeneous first integral of degree $d$ satisfies
\begin{equation}\label{eq:resonancemodp}
\alpha_{1}\rho_{1}+\dots+\alpha_{N}\rho_{N}\equiv d\pmod p .
\end{equation}
For instance, in characteristic $7$ the energy $\mathcal{E}$ admits a monomial with three factors from $\mathbb{A}_{e}(-1)$ and one from $\mathbb{A}_{e}(-\frac12)$, whose weight is $11\equiv4$. The reading of the Peirce data through resonances, which is the main point of Sections~\ref{sec:resonance} and \ref{sec:rac}, would then require qualifications without dynamical content. The algebraic theory over general fields, including the two rigidity statements of Subsection~\ref{ssec:1.2}, is developed in the forthcoming paper, which treats fields of characteristic not dividing $30$.
\end{remark}

\medskip\noindent\textbf{Principal notation.}~
For the convenience of the reader we list the main symbols and the places where they are introduced.

\begin{center}
\small
\begin{tabular}{p{0.25\textwidth}p{0.67\textwidth}}
\hline
$\mathbb{A}$, $\langle\cdot,\cdot\rangle$ & the algebra, its metric (invariant inner product); \eqref{eq:mult}\\
$L(x)$& the multiplication operator $L(x)y=xy$; \eqref{eq:Lx}\\
$q(x)$, $u(x)$, $\theta$ & the forms $\langle x,x\rangle$ and $\langle x,x^{2}\rangle$, and the constant in \eqref{eq:hsiang}\\
$\mathcal{R}$, $\mathcal{S}$ & the pseudocomposition and the admissible defects; \eqref{eq:defects}\\
$\operatorname{Idm}(\mathbb{A})$, $\operatorname{Nil}(\mathbb{A})$ & nonzero idempotents and $2$-nilpotents; Subsection~\ref{ssec:1.2}, Lemma~\ref{prop:nilpotent}\\
$\mathbb{A}_{e}(\mu)$, $n_{\mu}(e)$ & Peirce spaces and Peirce numbers of $e$; \eqref{eq:peirceno}\\
$\mathcal{K}_{e}$& the operator $\operatorname{id}-2L(e)$; \eqref{eq:Ke}\\
$\rho=1-2\mu$ & the Kovalevskaya exponents (eigenvalues of $\mathcal{K}_{e}$); \eqref{eq:kovalevskaya}\\
$\Sigma$, $\hat{u}$ & the unit sphere and the restriction of $u$ to it; Section~\ref{sec:morsebott}\\
$\delta(y)$, $\sigma(y)$ & the idempotent defect $y^{2}-y$ and the norm defect $\theta\langle y,y\rangle-1$; \eqref{eq:Dsigma}\\
$\mathbb{V}$, $\mathbb{V}^{\perp}$ & the Peirce space $\mathbb{A}_{e}(\frac12)$ and its orthogonal complement; \eqref{eq:VN}\\
$\mathcal{E}$, $\mathcal{F}$ & the energy and the sextic first integral of the flow $\dot x=x^{2}$; \eqref{eq:EFintro}\\
$\partial$ & the derivative along the flow, $(\partial f)(x)=df_{x}(x^{2})$; \eqref{eq:Lie}\\
\hline
\end{tabular}
\end{center}

\subsection{Organization}\label{ssec:1.6}
Section~\ref{sec:basic} collects the algebraic preliminaries, including the Peirce decomposition and two examples. Section~\ref{sec:euclid} proves Theorem~\ref{mthm:A}, and Section~\ref{sec:morsebott} proves Theorem~\ref{mthm:MB}. Section~\ref{sec:flow} proves Theorem~\ref{mthm:B}(i),(ii), and Section~\ref{sec:darboux} discusses Darboux points and proves Theorem~\ref{mthm:B}(iii). Section~\ref{sec:resonance} proves the resonance lemma and Theorem~\ref{mthm:C}, and Section~\ref{sec:rac} discusses the admissible side. Section~\ref{sec:final} contains concluding remarks and questions.

\medskip\noindent\textbf{Conventions.}~
Throughout, $\mathbb{A}$ is a finite-dimensional real algebra with a Euclidean metric (Remark~\ref{rem:standing}). If a polynomial function vanishes on $\mathbb{A}$, so do all its derivatives, and if it vanishes on a nonempty open subset, it vanishes identically. For a polynomial $f$ on $\mathbb{A}$ we write $df_{x}(y)$ for its differential at $x$ applied to $y$, and the gradient is defined by $\langle \nabla f(x),y\rangle =df_{x}(y)$. By invariance,
\begin{equation}\label{eq:dqdu}
dq_{x}(y)=2\langle x,y\rangle ,\qquad du_{x}(y)=3\langle x^{2},y\rangle ,
\end{equation}
that is,
\begin{equation}\label{eq:gradqu}
\nabla q(x)=2x,\qquad\nabla u(x)=3x^{2}.
\end{equation}

\section{Algebraic preliminaries}\label{sec:basic}

\subsection{The gradient identity}\label{ssec:2.1}
All the structure we shall use comes from one vector identity, obtained by differentiating \eqref{eq:hsiang}.

\begin{lemma}\label{lem:gradient}
Let $\mathbb{A}$ be a commutative algebra with a metric, let $\theta\in\mathbb{R}^{\times}$, and put
\begin{equation}\label{eq:P}
\mathcal{P}(x)=\langle x^{2},x^{3}\rangle -\theta q(x)u(x).
\end{equation}
Then
\begin{equation}\label{eq:gradP}
\nabla \mathcal{P}(x)=4\,x x^{3}+x^{2}x^{2}-2\theta u(x)\,x-3\theta q(x)\,x^{2}.
\end{equation}
In particular, if $\mathbb{A}$ is radial $\infty$-Hsiang, then
\begin{equation}\label{eq:V}
4\,x x^{3}+x^{2}x^{2}=2\theta u(x)\,x+3\theta q(x)\,x^{2}\qquad\text{for all }x\in\mathbb{A} .
\end{equation}
\end{lemma}
\begin{proof}
The differentials of the powers are $d(x^{2})_{x}(y)=2xy$ and $d(x^{3})_{x}(y)=2(xy)x+x^{2}y$. By invariance,
\begin{equation}\label{eq:invmoves}
\begin{aligned}
\langle xy,x^{3}\rangle &=\langle xx^{3},y\rangle ,\\
\langle x^{2},(xy)x\rangle &=\langle x^{3},xy\rangle =\langle xx^{3},y\rangle ,\\
\langle x^{2},x^{2}y\rangle &=\langle x^{2}x^{2},y\rangle .
\end{aligned}
\end{equation}
Hence
\begin{equation*}
d\langle x^{2},x^{3}\rangle _{x}(y)=2\langle xy,x^{3}\rangle +2\langle x^{2},(xy)x\rangle +\langle x^{2},x^{2}y\rangle =\langle 4xx^{3}+x^{2}x^{2},y\rangle .
\end{equation*}
By \eqref{eq:dqdu}, $d(qu)_{x}(y)=\langle 2u\,x+3q\,x^{2},y\rangle $. This gives \eqref{eq:gradP}. If $\mathcal{P}\equiv0$, then $\nabla \mathcal{P}\equiv0$, which is \eqref{eq:V}.
\end{proof}

\subsection{The Peirce decomposition}\label{ssec:2.2}
The following lemma is the reason why only the four values $1,-1,\pm\frac12$ occur in \eqref{eq:kovalevskaya}.

\begin{lemma}\label{prop:peirce}
Let $\mathbb{A}$ be radial $\infty$-Hsiang and $e\in\operatorname{Idm}(\mathbb{A})$. Then:
\begin{itemize}
\item[{(i)}] $\theta\langle e,e\rangle =1$; in particular $\mathbb{A}=\mathbb{R} e\oplus e^{\perp}$;
\item[{(ii)}] $e^{\perp}$ is invariant under $L(e)$, and
\begin{equation}\label{eq:peircepoly}
(L(e)+\operatorname{id})(2L(e)-\operatorname{id})(2L(e)+\operatorname{id})=0\qquad\text{on }e^{\perp};
\end{equation}
\item[{(iii)}] $L(e)$ is diagonalizable over $\mathbb{R}$, and
\begin{equation}\label{eq:peircedec}
\mathbb{A}=\mathbb{R} e\oplus\mathbb{A}_{e}(-1)\oplus\mathbb{A}_{e}(-\tfrac12)\oplus\mathbb{A}_{e}(\tfrac12)
\end{equation}
is an orthogonal decomposition, the metric being nondegenerate on each summand.
\end{itemize}
\end{lemma}
\begin{proof}
By Lemma~\ref{lem:gradient}, the right side of \eqref{eq:gradP} vanishes identically: $\nabla \mathcal{P}\equiv0$. Since $e^{k}=e$ for every $k\geq1$, we have $u(e)=q(e)=\langle e,e\rangle $, and therefore
\begin{equation*}
0=\nabla \mathcal{P}(e)=\bigl(4+1-2\theta\langle e,e\rangle -3\theta\langle e,e\rangle \bigr)e=5\bigl(1-\theta\langle e,e\rangle \bigr)e .
\end{equation*}
As $e\neq0$, this gives $\theta\langle e,e\rangle =1$. Since $\langle e,e\rangle \neq0$, we have $\mathbb{A}=\mathbb{R} e\oplus e^{\perp}$. For $y\perp e$, $\langle ey,e\rangle =\langle y,e^{2}\rangle =0$, so $e^{\perp}$ is $L(e)$-invariant.

Now differentiate $\nabla \mathcal{P}\equiv0$ at $e$ in a direction $y\in e^{\perp}$, and write $L=L(e)$. Using $\langle e,y\rangle =0$ we find
\begin{equation}\label{eq:dpowers}
\begin{aligned}
d(x^{3})_{e}(y)&=(2L^{2}+L)y,\\
d(xx^{3})_{e}(y)&=Ly+L(2L^{2}+L)y,\\
d(x^{2}x^{2})_{e}(y)&=4L^{2}y,\\
d(ux)_{e}(y)&=\langle e,e\rangle y,\\
d(qx^{2})_{e}(y)&=2\langle e,e\rangle Ly .
\end{aligned}
\end{equation}
Together with $\theta\langle e,e\rangle =1$ this gives
\begin{equation*}
0=d(\nabla \mathcal{P})_{e}(y)=(8L^{3}+8L^{2}+4L-2-6L)y=2(L+\operatorname{id})(2L-\operatorname{id})(2L+\operatorname{id})y .
\end{equation*}
The roots $-1,\frac12,-\frac12$ are distinct, so $L$ is diagonalizable on $e^{\perp}$ with eigenvalues among them. Eigenspaces of the self-adjoint operator $L$ for different eigenvalues are orthogonal, and since the metric is nondegenerate on $\mathbb{A}$, it is nondegenerate on each eigenspace.
\end{proof}

For $e\in\operatorname{Idm}(\mathbb{A})$ we put
\begin{equation}\label{eq:Ke}
\mathcal{K}_{e}=\operatorname{id}-2L(e).
\end{equation}
By Lemma~\ref{prop:peirce}, $\mathcal{K}_{e}$ is self-adjoint and acts on $\mathbb{A}_{e}(\mu)$ as the scalar $\rho=1-2\mu$: it equals $-1$ on $\mathbb{R} e$, and $3$, $2$, $0$ on $\mathbb{A}_{e}(-1)$, $\mathbb{A}_{e}(-\frac12)$, $\mathbb{A}_{e}(\frac12)$, respectively; in particular $\ker\mathcal{K}_{e}=\mathbb{A}_{e}(\frac12)$. This single operator will appear in three roles: $-\mathcal{K}_{e}$ is the differential at $e$ of the map $y\mapsto y^{2}-y$ (Section~\ref{sec:euclid}); $-3\sqrt{\theta}\,\mathcal{K}_{e}$ is the Hessian at $\sqrt{\theta}e$ of the cubic form restricted to the unit sphere (Section~\ref{sec:morsebott}); and $\mathcal{K}_{e}$ is the Kovalevskaya matrix at the Darboux point $-e$ of the flow \eqref{eq:flow}, with eigenvalues the Kovalevskaya exponents (Sections~\ref{sec:darboux} and \ref{sec:resonance}).

The eigenvalue $\frac12$ plays a special role in many classes of algebras with identities \cite{Tkachev21}. Here its meaning is geometric. The linearization at $e$ of the map $y\mapsto y^{2}-y$ is $-\mathcal{K}_{e}$, so $\mathbb{A}_{e}(\frac12)=\ker\mathcal{K}_{e}$ is the space of solutions of the linearized idempotent equation $2ey=y$, the natural candidate for the tangent space of $\operatorname{Idm}(\mathbb{A})$ at $e$; compare \cite{KrasnovTkachev19}. Theorem~\ref{mthm:A} shows that in the Euclidean case it is indeed the tangent space. As explained before Theorem~\ref{mthm:A}, this is far from automatic.

\subsection{Linearized defects}\label{ssec:2.3}
Next we compute the first-order behavior of the two defects $\mathcal{R}$ and $\mathcal{S}$ near an idempotent. This simple computation shows which Peirce eigenvalues are compatible with each of the vector identities \eqref{eq:radialclasses}.

\begin{lemma}\label{lem:defects}
Let $\mathbb{A}$ be radial $\infty$-Hsiang and $e\in\operatorname{Idm}(\mathbb{A})$. Then $\mathcal{R}(e)=\mathcal{S}(e)=0$, and for $y\in\mathbb{A}_{e}(\mu)$, $\mu\neq1$,
\begin{equation}\label{eq:dRdS}
\begin{aligned}
d\mathcal{R}_{e}(y)&=(2\mu-1)(\mu+1)\,y,\\
d\mathcal{S}_{e}(y)&=(2\mu-1)(2\mu+1)\,y .
\end{aligned}
\end{equation}
Consequently every $e\in\operatorname{Idm}(\mathbb{A})$ satisfies $n_{-1/2}(e)=0$ if $\mathbb{A}$ is radial pseudocomposition, and $n_{-1}(e)=0$ if $\mathbb{A}$ is radial admissible cubic.
\end{lemma}
\begin{proof}
By Lemma~\ref{prop:peirce}, $\mathcal{R}(e)=e-\theta\langle e,e\rangle e=0$ and $\mathcal{S}(e)=e-\theta u(e)e=0$. For $y\perp e$, by \eqref{eq:dqdu},
\begin{align*}
d\mathcal{R}_{e}(y)&=2(ey)e+ey-\theta\bigl(2\langle e,y\rangle e+\langle e,e\rangle y\bigr)=(2\mu^{2}+\mu-1)y,\\
d\mathcal{S}_{e}(y)&=2e(2ey)-\theta\bigl(3\langle e,y\rangle e+\langle e,e\rangle y\bigr)=(4\mu^{2}-1)y .
\end{align*}
If $\mathcal{R}\equiv0$, then $d\mathcal{R}_{e}=0$, which excludes $\mu=-\frac12$; if $\mathcal{S}\equiv0$, then $d\mathcal{S}_{e}=0$, which excludes $\mu=-1$.
\end{proof}

\subsection{Two examples}\label{ssec:2.4}
The following examples will be used to illustrate all main results. In both, the idempotent manifold has two components of different dimensions.

\begin{example}\label{ex:rpc}
Let $(\mathcal{H},h)$ be a Euclidean space of dimension $m\geq1$, and let $\mathbb{A}=\mathbb{R}\oplus\mathcal{H}$ with
\begin{equation}\label{eq:rpcmult}
\begin{aligned}
(r,x)(s,y)&=\bigl(rs-h(x,y),\,-(ry+sx)\bigr),\\
\langle (r,x),(s,y)\rangle &=rs+h(x,y).
\end{aligned}
\end{equation}
The metric is invariant, since both $\langle (r,x)(s,y),(p,z)\rangle $ and $\langle (r,x),(s,y)(p,z)\rangle $ are equal to
\begin{equation*}
rsp-r\,h(y,z)-s\,h(x,z)-p\,h(x,y).
\end{equation*}
For $z=(r,x)$ we have $z^{2}=(r^{2}-h(x,x),-2rx)$ and
\begin{equation*}
z^{3}=\bigl(r(r^{2}-h(x,x))+2r\,h(x,x),\;-(r^{2}-h(x,x))x+2r^{2}x\bigr)=\bigl(r^{2}+h(x,x)\bigr)(r,x),
\end{equation*}
so $\mathbb{A}$ is radial pseudocomposition with $\theta=1$. From $z^{2}=z$ we get $-2rx=x$ and $r^{2}-h(x,x)=r$. Hence
\begin{equation}\label{eq:rpcidm}
\operatorname{Idm}(\mathbb{A})=\{e_{0}\}\cup\{(-\tfrac12,v):h(v,v)=\tfrac34\},\qquad e_{0}=(1,0).
\end{equation}
Here $L(e_{0})(s,y)=(s,-y)$, so $(n_{-1},n_{-1/2},n_{1/2})(e_{0})=(m,0,0)$. For $e=(-\frac12,v)$ one computes
\begin{equation*}
L(e)(s,y)=\bigl(-\tfrac12s-h(v,y),\ \tfrac12y-sv\bigr),
\end{equation*}
so $L(e)$ equals $\frac12$ on $\{(0,y):h(v,y)=0\}$, and on the plane spanned by $(1,0)$ and $(0,v)$ it has the eigenvalues $1$ and $-1$. So $(n_{-1},n_{-1/2},n_{1/2})(e)=(1,0,m-1)$. Thus $\operatorname{Idm}(\mathbb{A})$ is the disjoint union of the isolated point $e_{0}$ and a sphere $S^{m-1}$ of dimension $m-1=n_{1/2}(e)$. The tangent space of this sphere at $(-\frac12,v)$ is $\{(0,y):h(v,y)=0\}=\mathbb{A}_{e}(\frac12)$, in agreement with Theorem~\ref{mthm:A}.
\end{example}

\begin{example}\label{ex:rac}
Let $\mathbb{A}$ be the space of symmetric $3\times3$ matrices over $\mathbb{R}$, with $\langle x,y\rangle =\operatorname{tr}(xy)$ and
\begin{equation}\label{eq:racmult}
x\cdot y=\tfrac12\bigl((x+y)^{\#}-x^{\#}-y^{\#}\bigr),
\end{equation}
where $x^{\#}$ is the adjugate matrix. Then $x\cdot x=x^{\#}$, and the metric is invariant because $\operatorname{tr}((x\cdot y)z)$ is the full polarization of the cubic form $3\det x$, hence symmetric in $x,y,z$. Moreover
\begin{equation}\label{eq:adjugate}
\begin{aligned}
(x^{\#})^{\#}&=\det(x)\,x,\\
\langle x,x^{\#}\rangle &=3\det x .
\end{aligned}
\end{equation}
Hence $\mathcal{S}\equiv0$ with $\theta=\frac13$, and $\mathbb{A}$ is radial admissible cubic; see \cite{NTV14,ElduqueOkubo00}. The nonzero idempotents are the symmetric matrices with $x^{\#}=x$. Diagonalizing $x=O\operatorname{diag}(a,b,c)O^{T}$, the equations $bc=a$, $ca=b$, $ab=c$ force $a,b,c\neq0$ (if one of them vanishes, all do); then $a^{2}=b^{2}=c^{2}=abc$, whence $a,b,c\in\{\pm1\}$ with $abc=1$. Thus
\begin{equation}\label{eq:racidm}
\operatorname{Idm}(\mathbb{A})=\{I\}\cup\{2ww^{T}-I:\ w^{T}w=1\}.
\end{equation}
One checks that $I\cdot y=\frac12(\operatorname{tr}(y)I-y)$, so $L(I)=-\frac12$ on trace-free matrices and the Peirce numbers $(n_{-1},n_{-1/2},n_{1/2})$ of $I$ are $(0,5,0)$, while $\operatorname{diag}(1,-1,-1)$ has Peirce numbers $(0,3,2)$. Thus $\operatorname{Idm}(\mathbb{A})$ is the union of a point and a real projective plane $\mathbb{R}P^{2}$ (the map $w\mapsto 2ww^{T}-I$ identifies $\pm w$), of dimensions $0$ and $2$.
\end{example}

\section{The idempotent manifold}\label{sec:euclid}

Positivity of the metric enters this section only through Lemma~\ref{prop:existence} (see Remark~\ref{rem:indefinite}). We write $|x|=\langle x,x\rangle ^{1/2}$.

\subsection{Existence and compactness}\label{ssec:3.1}
We begin by showing that a Euclidean radial $\infty$-Hsiang algebra always has idempotents and that they all lie on a sphere. The existence argument is variational. Idempotents arise from the maximum points of $u$ on the unit sphere.

\begin{lemma}\label{prop:existence}
Let $\mathbb{A}$ be a Euclidean radial $\infty$-Hsiang algebra with nonzero multiplication. Then $\operatorname{Idm}(\mathbb{A})\neq\emptyset$, $\theta>0$, and
\begin{equation}\label{eq:idmsphere}
\operatorname{Idm}(\mathbb{A})=\{x\in\mathbb{A}:\ x^{2}=x,\ \theta|x|^{2}=1\}
\end{equation}
is a compact real algebraic subset of the sphere of radius $\theta^{-1/2}$.
\end{lemma}
\begin{proof}
If $u\equiv0$, then the symmetric trilinear form $\langle xy,z\rangle $, which is a polarization of $u$, vanishes, and by nondegeneracy the multiplication is zero. So $u\not\equiv0$, and since $u$ is odd, its maximum $M$ on the unit sphere is positive. At a maximum point $x_{0}$, the Lagrange multiplier rule gives $3x_{0}^{2}=2\lambda x_{0}$ for some $\lambda\in\mathbb{R}$; pairing with $x_{0}$ gives $2\lambda=3M$, so $x_{0}^{2}=Mx_{0}$ and $e=M^{-1}x_{0}$ is a nonzero idempotent. By Lemma~\ref{prop:peirce}, $\theta|e|^{2}=1$, so $\theta>0$, and \eqref{eq:idmsphere} holds. The right side of \eqref{eq:idmsphere} is closed and bounded.
\end{proof}

\begin{remark}\label{rem:thetazero}
The existence argument does not use \eqref{eq:hsiang}: every commutative algebra with nonzero multiplication and a positive definite invariant inner product has a nonzero idempotent. This was observed by the second author in \cite[Lemma~2.1]{Tkachev16} and \cite{Tkachev18}; see also \cite{Tkachev19} for a more general setting. Consequently, if $u$ is a cubic form on a Euclidean space with $\Delta_{\infty}u\equiv0$, then $u\equiv0$. Indeed, otherwise the argument produces a nonzero idempotent $e$, while the computation in the proof of Lemma~\ref{prop:peirce} with $\theta=0$ gives $0=5e$. This is the cubic case of \cite[Proposition~4.1]{Tkachev16}, where it is shown, by a similar variational argument on the unit sphere, that a nonconstant real analytic $\infty$-harmonic function has no critical points. In the same paper, the existence of idempotents is the starting point of the proof that a cubic form solving the $p$-Laplace equation with $p\notin\{1,2\}$ vanishes identically \cite[Theorem~1.2]{Tkachev16}, which answers a question of J.~L.~Lewis in the cubic case. This is why $\theta\neq0$ is assumed throughout.
\end{remark}

\subsection{The idea of the proof}\label{ssec:3.2}
We look for the zeros of the polynomial map $\Phi(y)=y^{2}-y$ near a fixed idempotent $e$. By \eqref{eq:Ke}, its differential at $e$ is
\begin{equation}\label{eq:J}
d\Phi_{e}=2L(e)-\operatorname{id}=-\mathcal{K}_{e},
\end{equation}
a self-adjoint operator which vanishes on $\mathbb{V}=\mathbb{A}_{e}(\frac12)$ and is invertible on the orthogonal complement
\begin{equation}\label{eq:VN}
\mathbb{V}^{\perp}=\mathbb{R} e\oplus\mathbb{A}_{e}(-1)\oplus\mathbb{A}_{e}(-\tfrac12),
\end{equation}
where $\mathcal{K}_{e}$ has the eigenvalues $-1$, $3$, $2$ on the three summands. If $\mathbb{V}=0$, the inverse function theorem shows that $e$ is an isolated idempotent. If $\mathbb{V}\neq0$, the standard procedure is Lyapunov--Schmidt reduction (see, e.g., \cite[Section~2.4]{ChowHale82} for the general setting). In the finite-dimensional situation it amounts to an application of the implicit function theorem. One solves the $\mathbb{V}^{\perp}$-component of the equation $\Phi=0$ by the implicit function theorem, which expresses the $\mathbb{V}^{\perp}$-coordinate of a solution as a function $g(v)$ of its $\mathbb{V}$-coordinate. What remains is a system of $\dim\mathbb{V}$ equations in $\dim\mathbb{V}$ unknowns, the \emph{bifurcation equations}, which have no linear part. In general nothing more can be said. The bifurcation equations may have only the trivial solution, or a singular solution set, and the zero set of $\Phi$ near $e$ need not be a manifold.

\begin{figure}[ht]
\centering
\resizebox{0.95\textwidth}{!}{%
\begin{tikzpicture}[>=Stealth,font=\small]
\begin{scope}
\fill[blue!4,rounded corners=6pt] (-3.3,-2.3) rectangle (3.4,2.5);
\draw[->,black!60] (-3.0,0)--(3.1,0) node[above left,font=\small] {$\mathbb{V}=\mathbb{A}_{e}(\frac12)$};
\draw[->,black!60] (0,-2.0)--(0,2.2) node[below right,font=\small] {$\mathbb{V}^{\perp}$};
\draw[very thick,blue!60!black,domain=-2.7:2.7,samples=60] plot (\x,{-0.17*\x*\x});
\node[blue!60!black,font=\small,anchor=north] at (-1.9,-0.75) {$w=g(v)$};
\fill (0,0) circle (2pt) node[above left,font=\small] {$e$};
\coordinate (Y) at (1.8,{-0.17*1.8*1.8});
\fill[blue!60!black] (Y) circle (2pt);
\node[font=\small,anchor=north] at (1.75,-0.62) {$y(v)$};
\draw[->,thick,orange!80!black] (Y) -- ++(0.95,0) node[below right=-2pt,font=\small,text=orange!60!black] {$\delta\in\mathbb{V}$};
\draw[black!40,dashed] (1.8,0)--(Y);
\node[font=\small] at (1.8,-0.2) {};
\fill[black!50] (-1.2,1.1) circle (1.6pt);
\node[font=\scriptsize,align=left,anchor=west,text=black!60] at (-1.05,1.25) {off the graph:\\ $\pi_{\mathbb{V}^{\perp}}\delta\neq0$};
\node[font=\small\bfseries] at (0,2.85) {(a) the reduction};
\end{scope}
\begin{scope}[xshift=7.6cm,yshift=-1.9cm]
\fill[orange!5,rounded corners=6pt] (-0.7,-0.4) rectangle (6.4,4.4);
\node[anchor=west,font=\small] (r1) at (-0.45,3.55) {$\deg 1:\ \ 6\delta_{1}-5\sigma_{1}e\;=\;0$};
\node[anchor=west,font=\small] (r2) at (-0.45,2.6) {$\deg 2:\ \ 6\delta_{2}-5\sigma_{2}e\;=\;Q_{2}(\delta_{1},\sigma_{1})$};
\node[anchor=west,font=\small] at (1.0,1.95) {$\vdots$};
\node[anchor=west,font=\small] (rk) at (-0.45,1.25) {$\deg k:\ \ 6\delta_{k}-5\sigma_{k}e\;=\;Q_{k}(\delta_{j},\sigma_{j};\,j<k)$};
\draw[->,thick,green!45!black] (3.05,3.35) to[bend left=25] (4.1,2.8);
\draw[->,thick,green!45!black] (4.75,2.4) to[bend left=25] (5.2,1.5);
\node[font=\scriptsize,text=green!35!black,anchor=west] at (4.15,3.25) {$=0$};
\node[font=\scriptsize,align=center,text=black!75] at (2.85,0.25) {$6\delta_{k}\in\mathbb{V},\ \ 5\sigma_{k}e\in\mathbb{R}e\subset\mathbb{V}^{\perp}$\\[1pt] $\Longrightarrow\ \delta_{k}=0,\ \sigma_{k}=0$};
\node[font=\small\bfseries] at (2.85,4.75) {(b) the defect identity};
\end{scope}
\end{tikzpicture}
}
\caption{The two ingredients of the proof of Theorem~\ref{mthm:A} near an idempotent $e$. (a) The Lyapunov--Schmidt reduction (Lemma~\ref{lem:reduction}): the $\mathbb{V}^{\perp}$-component of the equation $y^{2}=y$ is solved by the implicit function theorem, which gives the graph $w=g(v)$ over $\mathbb{V}=\mathbb{A}_{e}(\frac12)$; at the points $y(v)$ of this graph the idempotent defect $\delta=y^{2}-y$ lies in $\mathbb{V}$. (b) The defect identity (Lemma~\ref{lem:defectidentity}) along the graph, compared degree by degree in $v$: in degree $k$ its linear part is $6\delta_{k}-5\sigma_{k}e$, while all remaining terms, collected in $Q_{k}$, involve only the defects of lower degree. Since $6\delta_{k}\in\mathbb{V}$ and $5\sigma_{k}e\perp\mathbb{V}$, induction on $k$ shows that $\sigma$ and $\delta$ vanish to infinite order at $e$; being real analytic, they vanish identically, so the whole graph consists of idempotents.}
\label{fig:proofA}
\end{figure}
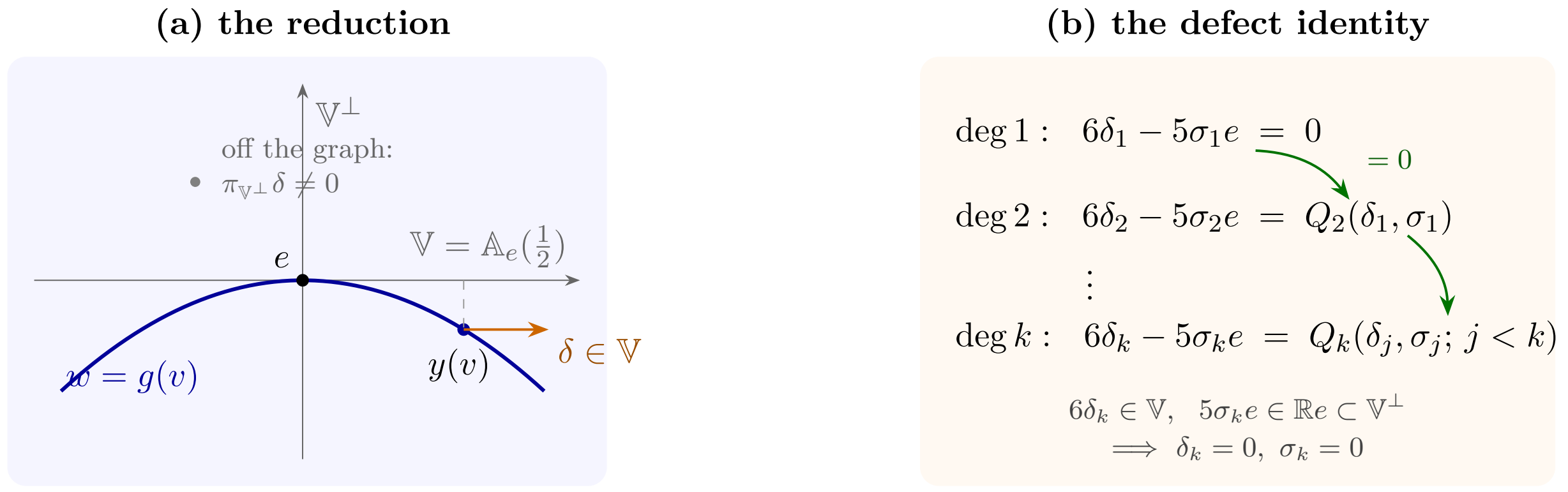

In a radial $\infty$-Hsiang algebra the bifurcation equations hold automatically. The reason is an identity (Lemma~\ref{lem:defectidentity}) which expresses a certain linear combination of the idempotent defect $\delta=y^{2}-y$ and of the norm defect $\sigma=\theta|y|^{2}-1$ through terms quadratic in $(\sigma,\delta)$. At $y=e$ this linear combination is nondegenerate precisely in the $\mathbb{V}$-directions. Along the graph $y(v)=e+v+g(v)$ provided by the reduction, the idempotent defect $\delta(y(v))$ lies in $\mathbb{V}$. We expand $\delta(y(v))$ and $\sigma(y(v))$ in homogeneous parts $\delta_{k}$, $\sigma_{k}$ in $v$ and compare the parts of degree $k$ on both sides of the defect identity. On the left, this part equals $6\delta_{k}-5\sigma_{k}e$ plus terms that contain some $\delta_{j}$ or $\sigma_{j}$ with $j<k$; on the right, it is a sum of products of parts of lower degree. Since $\delta_{k}\in\mathbb{V}$ and $e\perp\mathbb{V}$, the map $(\sigma_{k},\delta_{k})\mapsto6\delta_{k}-5\sigma_{k}e$ is injective on $\mathbb{R}\times\mathbb{V}$, so if all parts of degree less than $k$ vanish, then so do $\delta_{k}$ and $\sigma_{k}$.
By induction on the degree all of them vanish, and since $\delta$ and $\sigma$ are real analytic, they vanish identically. The zero set of $\Phi$ near $e$ is then the graph of $g$, which is a manifold with tangent space $\mathbb{V}$. Figure~\ref{fig:proofA} illustrates the two steps.

\subsection{The reduction}\label{ssec:3.3}
Fix $e\in\operatorname{Idm}(\mathbb{A})$, let $\mathbb{V}=\mathbb{A}_{e}(\frac12)$ and let $\mathbb{V}^{\perp}$ be its orthogonal complement; since $L(e)$ is self-adjoint, $\mathbb{V}^{\perp}$ is the range of $\mathcal{K}_{e}$, on which $\mathcal{K}_{e}$ is invertible, and in a radial $\infty$-Hsiang algebra $\mathbb{V}^{\perp}$ is given by \eqref{eq:VN}. Let  $\pi_{\mathbb{V}}$ and $\pi_{\mathbb{V}^{\perp}}$ be the orthogonal projections onto $\mathbb{V}$ and $\mathbb{V}^{\perp}$. They commute with $L(e)$. For $y\in\mathbb{A}$ put
\begin{equation}\label{eq:Dsigma}
\begin{aligned}
\delta(y)&=y^{2}-y,\\
\sigma(y)&=\theta\langle y,y\rangle -1 ,
\end{aligned}
\end{equation}
the \emph{idempotent defect} and the \emph{norm defect} of $y$. Write $y=e+v+w$ with $v=\pi_{\mathbb{V}}(y-e)=\pi_{\mathbb{V}}(y)\in\mathbb{V}$ and $w=\pi_{\mathbb{V}^{\perp}}(y-e)\in\mathbb{V}^{\perp}$. Using $e^{2}=e$ and $2ev=v$ yields
\begin{equation}\label{eq:split}
\delta(e+v+w)=e+2ev+2ew+(v+w)^{2}-e-v-w=-\mathcal{K}_{e}w+(v+w)^{2}.
\end{equation}
The term $\mathcal{K}_{e}w$ lies in $\mathbb{V}^{\perp}$, so $\delta(e+v+w)=0$ if and only if
\begin{equation}\label{eq:system}
\begin{aligned}
\mathcal{K}_{e}w-\pi_{\mathbb{V}^{\perp}}\bigl((v+w)^{2}\bigr)&=0,\\
\pi_{\mathbb{V}}\bigl((v+w)^{2}\bigr)&=0 .
\end{aligned}
\end{equation}
The first line in \eqref{eq:system} is uniquely solvable for $w$ near $0$, because its linearization in $w$ is the invertible operator $\mathcal{K}_{e}|_{\mathbb{V}^{\perp}}$. The second one is the bifurcation equation.

\begin{lemma}\label{lem:reduction}
Let $\mathbb{A}$ be a commutative algebra with a metric, and let $e\in\operatorname{Idm}(\mathbb{A})$. There are open balls $U\subset\mathbb{V}$ and $U'\subset\mathbb{V}^{\perp}$ centered at $0$ and a real analytic map $g:U\to U'$ with $g(0)=0$ and $dg(0)=0$ with the following properties.
\begin{itemize}
\item[{(a)}] For $(v,w)\in U\times U'$, the first line of \eqref{eq:system} holds if and only if $w=g(v)$.
\item[{(b)}] For $v\in U$ put $y(v)=e+v+g(v)$. Then the idempotent defect of $y(v)$ lies in $\mathbb{V}$:
\begin{equation}\label{eq:DinV}
\delta(y(v))=\pi_{\mathbb{V}}\bigl((v+g(v))^{2}\bigr)\in\mathbb{V}.
\end{equation}
\item[{(c)}] Near $e$, the idempotents are exactly the points $y(v)$ with vanishing idempotent defect:
\begin{equation}\label{eq:localgraph}
\operatorname{Idm}(\mathbb{A})\cap(e+U+U')=\{y(v):\ v\in U,\ \delta(y(v))=0\}.
\end{equation}
\end{itemize}
\end{lemma}
\begin{proof}
Consider the polynomial map
\begin{equation*}
G:\mathbb{V}\times\mathbb{V}^{\perp}\to\mathbb{V}^{\perp},\qquad G(v,w)=\mathcal{K}_{e}w-\pi_{\mathbb{V}^{\perp}}\bigl((v+w)^{2}\bigr).
\end{equation*}
We have $G(0,0)=0$, and the partial differential $G_{w}(0,0)=\mathcal{K}_{e}|_{\mathbb{V}^{\perp}}$ is invertible. The real analytic implicit function theorem (see, e.g., \cite{KrantzParks02}) gives balls $U$, $U'$ and a real analytic $g$ such that (a) holds. Since $G(v,0)=-\pi_{\mathbb{V}^{\perp}}(v^{2})$ is quadratic in $v$, we have $G_{v}(0,0)=0$, and differentiating $G(v,g(v))=0$ at $v=0$ gives $\mathcal{K}_{e}\,dg(0)=0$, hence $dg(0)=0$. By \eqref{eq:split} and (a),
\begin{equation*}
\delta(y(v))=-\mathcal{K}_{e}g(v)+(v+g(v))^{2}=-G(v,g(v))+\pi_{\mathbb{V}}\bigl((v+g(v))^{2}\bigr)=\pi_{\mathbb{V}}\bigl((v+g(v))^{2}\bigr),
\end{equation*}
which is (b). Finally, for $(v,w)\in U\times U'$, the point $e+v+w$ is idempotent if and only if both lines of \eqref{eq:system} hold, that is, if and only if $w=g(v)$ and $\delta(y(v))=0$. Shrinking the balls we may assume $0\notin e+U+U'$, so that all these idempotents are nonzero. This proves (c).
\end{proof}

\subsection{The defect identity}\label{ssec:3.4}
The next identity holds for \emph{all} $y$, not only near $e$. It is simply the gradient identity \eqref{eq:V} rewritten in terms of the two defects \eqref{eq:Dsigma}.

\begin{lemma}[Defect identity]\label{lem:defectidentity}
Let $\mathbb{A}$ be radial $\infty$-Hsiang. For every $y\in\mathbb{A}$, with $\delta=\delta(y)$ and $\sigma=\sigma(y)$,
\begin{equation}\label{eq:defectidentity}
2\delta+6\,y\delta+4\,y(y\delta)-2\theta\langle y,\delta\rangle \,y-5\sigma\,y=3\sigma \delta-\delta^{2}.
\end{equation}
\end{lemma}
\begin{proof}
Substitute $y^{2}=y+\delta$ into \eqref{eq:V}. We have
\begin{equation}\label{eq:powersD}
\begin{aligned}
y^{3}&=y+\delta+y\delta,\\
yy^{3}&=y+\delta+y\delta+y(y\delta),\\
y^{2}y^{2}&=y+\delta+2y\delta+\delta^{2},
\end{aligned}
\end{equation}
and
\begin{equation}\label{eq:scalarsD}
\begin{aligned}
\theta q(y)&=1+\sigma,\\
\theta u(y)&=\theta\langle y,y+\delta\rangle =1+\sigma+\theta\langle y,\delta\rangle .
\end{aligned}
\end{equation}
Insert \eqref{eq:powersD} and \eqref{eq:scalarsD} into $4yy^{3}+y^{2}y^{2}-2\theta u(y)\,y-3\theta q(y)\,y^{2}=0$. The terms involving neither $\delta$ nor $\sigma$ give $(4+1-2-3)y=0$. The remaining multiples of $y$ are $-5\sigma y-2\theta\langle y,\delta\rangle y$, and the other terms are
\begin{equation*}
4\delta+4y\delta+4y(y\delta)+\delta+2y\delta+\delta^{2}-3(1+\sigma)\delta=2\delta+6y\delta+4y(y\delta)+\delta^{2}-3\sigma \delta .
\end{equation*}
Collecting all terms gives \eqref{eq:defectidentity}.
\end{proof}

It is instructive to look at the left side of \eqref{eq:defectidentity} at $y=e$, as a linear function of $(\sigma,\delta)$. On $\mathbb{A}_{e}(\mu)$ the operator
\begin{equation}\label{eq:complementarity}
2\operatorname{id}+6L(e)+4L(e)^{2}=2\bigl(2L(e)+\operatorname{id}\bigr)\bigl(L(e)+\operatorname{id}\bigr)
\end{equation}
acts by the scalar $2(2\mu+1)(\mu+1)$. It vanishes on $\mathbb{A}_{e}(-1)\oplus\mathbb{A}_{e}(-\frac12)$, where the equation $y^{2}=y$ is nondegenerate, and it equals $6$ on $\mathbb{V}=\mathbb{A}_{e}(\frac12)$, where the equation is degenerate. In particular, for $s\in\mathbb{R}$ and $\eta\in\mathbb{V}$ the left side of \eqref{eq:defectidentity} with $y=e$, $\sigma=s$, $\delta=\eta$ equals
\begin{equation}\label{eq:K0}
\bigl(2+6\cdot\tfrac12+4\cdot\tfrac14\bigr)\eta-5se=6\eta-5se ,
\end{equation}
since $\langle e,\eta\rangle =0$.
This complementarity between the degenerate directions of $\mathcal{K}_{e}$ and the nondegenerate directions of \eqref{eq:complementarity} is the heart of the proof.

\subsection{Proof of Theorem~\ref{mthm:A}}\label{ssec:3.5}
We now combine the reduction lemma with the defect identity. The idea is to compare the homogeneous parts of degree $k$ in $v$ of the two sides of \eqref{eq:defectidentity} at the point $y(v)$. The left side is linear in the defects $(\sigma,\delta)$, and its part of degree $k$ is the injective map \eqref{eq:K0} applied to the parts of degree $k$ of the defects, up to terms of lower degree; the right side is quadratic and involves only parts of lower degree. Hence, by induction on the degree, both defects vanish to infinite order at $v=0$, and therefore identically, since they are real analytic. This gives the following local statement, from which Theorem~\ref{mthm:A} follows at once.

\begin{theorem}\label{thm:idmsmooth}
Let $\mathbb{A}$ be a Euclidean radial $\infty$-Hsiang algebra with nonzero multiplication, and let $e\in\operatorname{Idm}(\mathbb{A})$. With the notation of Lemma~\ref{lem:reduction}, $\delta(y(v))=0$ and $\sigma(y(v))=0$ for all $v\in U$. Consequently, near $e$, $\operatorname{Idm}(\mathbb{A})$ is the graph $\{e+v+g(v):v\in U\}$ of a real analytic map; it is a real analytic submanifold of dimension $n_{1/2}(e)$, and its tangent space at $e$ is $\mathbb{A}_{e}(\frac12)$.
\end{theorem}
\begin{proof}
For $v\in U$ write $y(v)=e+v+g(v)$, $\delta(v)=\delta(y(v))$ and $\sigma(v)=\sigma(y(v))$. Since $g$ is real analytic, so are $\delta$ and $\sigma$. First we expand them at $v=0$ in homogeneous polynomials:
\begin{equation}\label{eq:expansions}
y(v)=e+\sum_{k\geq1}y_{k}(v),\qquad
\delta(v)=\sum_{k\geq0}\delta_{k}(v),\qquad
\sigma(v)=\sum_{k\geq0}\sigma_{k}(v),
\end{equation}
where $y_{1}(v)=v$ and $y_{k}$ is the part of degree $k$ of $g$ for $k\geq2$. The part of degree $k$ of a product (in $\mathbb{A}$, of the inner product, or with a scalar function) of two such series $\sum_{i}a_{i}$ and $\sum_{j}b_{j}$ is $\sum_{i+j=k}a_{i}b_{j}$. By \eqref{eq:DinV}, $\delta(v)=\pi_{\mathbb{V}}\bigl((v+g(v))^{2}\bigr)$ with $v+g(v)=\sum_{i\geq1}y_{i}(v)$, so
\begin{equation}\label{eq:deltak}
\delta_{0}=\delta_{1}=0,\qquad \delta_{k}=\pi_{\mathbb{V}}\Bigl(\sum_{i+j=k,\ i,j\geq1}y_{i}y_{j}\Bigr)\in\mathbb{V}\quad(k\geq2).
\end{equation}
Moreover $\sigma_{0}=\sigma(e)=0$ by Lemma~\ref{prop:peirce}.

Next we show, by induction on $k$, that $\delta_{k}=0$ and $\sigma_{k}=0$ for all $k\geq0$. For $k=0$ this was just noted. Let $k\geq1$ and suppose that $\delta_{j}=0$ and $\sigma_{j}=0$ for $j<k$. We compare the parts of degree $k$ of the two sides of \eqref{eq:defectidentity} at $y=y(v)$. Every term of degree $k$ on the right side $3\sigma\delta-\delta^{2}$ has the form $3\sigma_{i}\delta_{j}$ or $\delta_{i}\delta_{j}$ with $i+j=k$; as $\sigma_{0}=\delta_{0}=0$, it contains a factor of degree less than $k$, so it vanishes. On the left side, each term is linear in $(\sigma,\delta)$; replacing $y$ by $e+\sum_{i\geq1}y_{i}$, its part of degree $k$ is the same expression with $y$ replaced by $e$ and $(\sigma,\delta)$ by $(\sigma_{k},\delta_{k})$, plus terms containing some $\sigma_{j}$ or $\delta_{j}$ with $j<k$, which vanish. Since $\delta_{k}\in\mathbb{V}$, \eqref{eq:K0} gives
\begin{equation}\label{eq:degk}
6\delta_{k}-5\sigma_{k}e=0 .
\end{equation}
Here $\delta_{k}\in\mathbb{V}$ and $\sigma_{k}e\in\mathbb{R}e\subset\mathbb{V}^{\perp}$, so $\delta_{k}=0$ and $\sigma_{k}=0$, which completes the induction.

Thus all Taylor coefficients at $0$ of the real analytic maps $\delta$ and $\sigma$ vanish. Since $U$ is connected, $\delta\equiv0$ and $\sigma\equiv0$ on $U$ (see \cite{KrantzParks02}). By \eqref{eq:localgraph}, $\operatorname{Idm}(\mathbb{A})\cap(e+U+U')=\{y(v):v\in U\}$ is the graph of the real analytic map $g$ over the ball $U\subset\mathbb{V}$. A graph of a real analytic map is a real analytic submanifold of dimension $\dim\mathbb{V}=n_{1/2}(e)$, and its tangent space at $e$ is $\mathbb{V}$ because $dg(0)=0$.
\end{proof}

\begin{proof}[Proof of Theorem~\ref{mthm:A}]
By Lemma~\ref{prop:existence}, $\operatorname{Idm}(\mathbb{A})$ is nonempty and compact and lies on the sphere \eqref{eq:sphere}; by Theorem~\ref{thm:idmsmooth} it is a real analytic submanifold whose tangent space at $e$ is $\mathbb{A}_{e}(\frac12)$. Hence $n_{1/2}$ is the local dimension, and it is constant on connected components. For the other two Peirce numbers we use the trace. The function $e'\mapsto\operatorname{tr} L(e')$ is continuous on $\operatorname{Idm}(\mathbb{A})$, and by \eqref{eq:peircedec}
\begin{equation}\label{eq:trace}
\operatorname{tr} L(e')=1-n_{-1}(e')-\tfrac12n_{-1/2}(e')+\tfrac12n_{1/2}(e')\in\tfrac12\mathbb{Z}.
\end{equation}
A continuous function with values in the discrete set $\frac12\mathbb{Z}$ is constant on each connected component. Since also
\begin{equation*}
n_{-1}(e')+n_{-1/2}(e')=\dim\mathbb{A}-1-n_{1/2}(e'),
\end{equation*}
the numbers $n_{-1}(e')$ and $n_{-1/2}(e')$ are constant on each component as well.
\end{proof}

\begin{remark}
The sphere constraint played no role in the construction of $y(v)$: we obtained $\sigma(y(v))=0$ as an output of the defect identity, at the same time as $\delta(y(v))=0$. The argument uses only the vector identity \eqref{eq:V} and Lemma~\ref{prop:peirce}.
\end{remark}

\begin{remark}\label{rem:indefinite}
Indefinite metrics lie outside the scope of this paper, but the local part of the argument does not depend on positivity. Positivity of the metric was used only in Lemma~\ref{prop:existence}. In Lemma~\ref{lem:reduction} and Theorem~\ref{thm:idmsmooth} one only needs the projections $\pi_{\mathbb{V}}$, $\pi_{\mathbb{V}^{\perp}}$, which exist because the metric is nondegenerate on each Peirce space, and the fact that $e\notin\mathbb{V}$. Hence, for any radial $\infty$-Hsiang algebra over $\mathbb{R}$ with a nondegenerate (possibly indefinite) metric, $\operatorname{Idm}(\mathbb{A})$ is near each of its points a real analytic manifold of dimension $n_{1/2}(e)$ with tangent space $\mathbb{A}_{e}(\frac12)$; over $\mathbb{C}$, by the holomorphic implicit function theorem, it is a complex manifold with the same tangent spaces. This is not the case for Lemma~\ref{prop:nilpotent} below, whose proof uses positivity.
\end{remark}

\begin{remark}[Power series description]
The map $g$ can be computed recursively. Write $y=e+v+\sum_{k\geq2}g_{k}(v)$ with $g_{k}$ homogeneous of degree $k$ with values in $\mathbb{V}^{\perp}$, and put $g_{1}(v)=v$. Comparing the terms of degree $k\geq2$ in $y^{2}=y$ gives, with sums over $i,j\geq1$,
\begin{equation}\label{eq:recursion}
\begin{aligned}
&\mathcal{K}_{e}\,g_{k}=\pi_{\mathbb{V}^{\perp}}\Bigl(\sum_{i+j=k}g_{i}g_{j}\Bigr),\\
&\pi_{\mathbb{V}}\Bigl(\sum_{i+j=k}g_{i}g_{j}\Bigr)=0.
\end{aligned}
\end{equation}
The first equation determines $g_{k}$; the second is a compatibility condition which, by Theorem~\ref{thm:idmsmooth}, holds automatically in every degree. For instance $g_{2}=\mathcal{K}_{e}^{-1}\pi_{\mathbb{V}^{\perp}}(v^{2})$, and the condition in degree $2$ says that $\pi_{\mathbb{V}}(v^{2})=0$, that is, $\mathbb{A}_{e}(\frac12)^{2}\subset\mathbb{A}_{e}(\frac12)^{\perp}$.
\end{remark}

\begin{remark}
In Example~\ref{ex:rpc}, the component $\{(-\frac12,v)\}$ is a sphere of dimension $m-1=n_{1/2}$, and the isolated point $e_{0}$ has $n_{1/2}(e_{0})=0$. In Example~\ref{ex:rac}, the components are the point $\{I\}$ with $n_{1/2}(I)=0$ and a projective plane with $n_{1/2}=2$. In both examples the Peirce numbers differ between components, so ``constant on components'' cannot be improved to ``constant'' in general. For exact Hsiang algebras the Peirce numbers are known to be globally constant (this is again proved in the monograph of the second author mentioned above and is not used here).
\end{remark}

\begin{remark}\label{rem:counter}
The identity \eqref{eq:hsiang} cannot be dropped in Theorem~\ref{thm:idmsmooth}, even if $\operatorname{Idm}(\mathbb{A})$ is a real analytic manifold. Let $\mathbb{A}=\mathbb{R}\oplus\mathbb{R}^{3}$ with the inner product of Example~\ref{ex:rpc} ($m=3$, $h$ Euclidean), write $x=(x_{1},x_{2},x_{3})$, $y=(y_{1},y_{2},y_{3})$, identify $(x_{2},x_{3})$ with $z=x_{2}+ix_{3}\in\mathbb{C}$ and $(y_{2},y_{3})$ with $w=y_{2}+iy_{3}$, and perturb \eqref{eq:rpcmult} to
\begin{equation*}
(r,x)(s,y)=\bigl(rs-\langle x,y\rangle,\,-(ry+sx)\bigr)+\bigl(0,0,\overline{zw}\bigr).
\end{equation*}
The metric is still invariant, because $\langle\overline{zw},\zeta\rangle=\operatorname{Re}(zw\zeta)$ is symmetric in $z,w,\zeta$, and the cubic form is $u=r^{3}-3r|x|^{2}+\operatorname{Re}z^{3}$. The equation $(r,x)^{2}=(r,x)$ reads
\begin{equation*}
r^{2}-|x|^{2}=r,\qquad x_{1}(1+2r)=0,\qquad \bar z^{\,2}=(1+2r)z .
\end{equation*}
If $z=0$, we get $e_{0}=(1,0)$ and $e_{\pm}=(-\frac12,\pm\frac{\sqrt3}{2},0,0)$. If $z\neq0$, then $|z|=|1+2r|\neq0$, so $x_{1}=0$, $z=(1+2r)\omega$ with $\omega^{3}=1$, and $3r^{2}+5r+1=0$. Hence $\operatorname{Idm}(\mathbb{A})$ consists of $9$ points and is a $0$-dimensional analytic manifold. On the other hand, the perturbation vanishes when $z=0$, so $L(e_{\pm})$ is the same as in Example~\ref{ex:rpc}: $\mathbb{V}=\mathbb{A}_{e_{\pm}}(\frac12)$ is the $z$-plane and $n_{1/2}(e_{\pm})=2\neq0=\dim_{e_{\pm}}\operatorname{Idm}(\mathbb{A})$. In terms of Lemma~\ref{lem:reduction}, which does not use \eqref{eq:hsiang}, the bifurcation equation \eqref{eq:DinV} at $e_{\pm}$ becomes $\bar z^{\,2}=0$, its leading term being $\pi_{\mathbb{V}}(v^{2})=\bar z^{\,2}=\frac13\nabla\operatorname{Re}z^{3}$; this is exactly the degree-$2$ compatibility condition in \eqref{eq:recursion} that fails. Accordingly, $\mathbb{A}$ is not radial $\infty$-Hsiang: $\langle x^{2},x^{3}\rangle=\theta\langle x,x\rangle\langle x,x^{2}\rangle$ forces $\theta=1$ at $x=(1,0)$, but $\theta=3$ at $x=(0;0,1,0)$, where $x^{2}=(-1;0,1,0)$ and $x^{3}=(-1;0,2,0)$.
\end{remark}

\section{The cubic form on the sphere: a Morse--Bott picture}\label{sec:morsebott}

For any commutative algebra with an invariant inner product, the critical points on the unit sphere of its cubic form $u$ are rescaled nonzero idempotents and their negatives and nonzero $2$-nilpotents (Lemma~\ref{lem:critical} below; compare \cite{KrasnovTkachev19,NTV14}). For a general cubic form, the restriction of $u$ to the sphere is a complicated function, with critical points having different indices and with possibly degenerate critical sets. A classical situation where the picture is rigid is that of isoparametric hypersurfaces. By M\"unzner's theorem \cite{Munzner80}, the restriction to the sphere of a Cartan--M\"unzner polynomial is a Morse--Bott function whose only critical sets are the two focal submanifolds, on which it attains its maximum and its minimum. The simplest instances are the Cartan isoparametric cubics \cite{Cartan39}. They satisfy $|\nabla u|^{2}=9|x|^{4}$, that is, $\mathcal{E}\equiv0$ with $\theta=1$ in the notation of \eqref{eq:EFintro}, so by Theorem~\ref{thm:integrals}(iii) their algebras are radial pseudocomposition. They have no nonzero $2$-nilpotents, and, up to scaling, their idempotent manifolds are the focal submanifolds on which $u$ attains its maximum, which are the Veronese embeddings of the projective planes over $\mathbb{R}$, $\mathbb{C}$, $\mathbb{H}$ and $\mathbb{O}$; see \cite{Cartan39,NTV14}. Here $\dim\mathbb{A}=3k+2$ with $k=1,2,4,8$, and every idempotent has Peirce numbers $(n_{-1},n_{-1/2},n_{1/2})=(k+1,0,2k)$, in accord with Theorem~\ref{mthm:A}. 

Theorem~\ref{mthm:MB} shows that any radial $\infty$-Hsiang algebra exhibits a similar rigidity: $u$ has no positive critical points on the sphere other than its maximum points, and its restriction to the sphere is a Morse--Bott function, whose critical set consists of the maximum set, the minimum set and normalized nonzero $2$-nilpotents. The new feature is that the maximum set may have several components of different dimensions, as in Examples~\ref{ex:rpc} and \ref{ex:rac}.

In this section $\mathbb{A}$ is a Euclidean radial $\infty$-Hsiang algebra with nonzero multiplication, $\Sigma=\{x:|x|=1\}$ is the unit sphere and $\hat{u}=u|_{\Sigma}$. By \eqref{eq:gradqu}, the flow \eqref{eq:flow} is the gradient flow of $\frac13u$, and its radial projection to $\Sigma$ is, up to a change of time, the gradient flow of $\hat{u}$. This links the geometry of this section with the dynamics described in Sections~\ref{sec:flow}--\ref{sec:resonance}.

\begin{lemma}\label{lem:critical}
A point $x\in \Sigma$ is a critical point of $\hat{u}$ if and only if $x^{2}=u(x)\,x$. In this case:
\begin{itemize}
\item[{(i)}] if $u(x)>0$, then $x=\sqrt{\theta}\,e$ for a unique $e\in\operatorname{Idm}(\mathbb{A})$, and $u(x)=\sqrt{\theta}$;
\item[{(ii)}] if $u(x)=0$, then $x^{2}=0$;
\item[{(iii)}] if $u(x)<0$, then $-x$ is a critical point with $u(-x)>0$.
\end{itemize}
Conversely, $\sqrt{\theta}\,e\in \Sigma$ is a critical point of $\hat{u}$ for every $e\in\operatorname{Idm}(\mathbb{A})$.
\end{lemma}
\begin{proof}
By the Lagrange multiplier rule, $x\in \Sigma$ is critical if and only if $3x^{2}=2\lambda x$ for some $\lambda\in\mathbb{R}$; pairing with $x$ gives $2\lambda=3u(x)$, so the condition is $x^{2}=u(x)x$. If $u(x)>0$, put $e=u(x)^{-1}x$. Then $e^{2}=u(x)^{-2}x^{2}=e$, and by Lemma~\ref{prop:peirce}, $1=\theta|e|^{2}=\theta u(x)^{-2}$, so $u(x)=\sqrt{\theta}$ and $x=\sqrt{\theta}\,e$. The case $u(x)=0$ is immediate, and the case $u(x)<0$ follows because $u$ is odd. Conversely, if $e\in\operatorname{Idm}(\mathbb{A})$ and $x=\sqrt{\theta}e$, then $|x|=1$, $x^{2}=\theta e=\sqrt{\theta}\,x$ and $u(x)=\theta^{3/2}\langle e,e\rangle =\sqrt{\theta}$.
\end{proof}

\begin{lemma}\label{prop:hessianU}
Let $e\in\operatorname{Idm}(\mathbb{A})$ and $x_{0}=\sqrt{\theta}\,e$. The Hessian of $\hat{u}$ at $x_{0}$, as a quadratic form on $T_{x_{0}}\Sigma=e^{\perp}$, is
\begin{equation}\label{eq:hessU}
\operatorname{Hess}_{x_{0}}\hat{u}(y,y)=-3\sqrt{\theta}\,\langle \mathcal{K}_{e}y,y\rangle ,\qquad y\in e^{\perp}.
\end{equation}
It vanishes on $\mathbb{A}_{e}(\frac12)$, equals $-6\sqrt{\theta}\,|y|^{2}$ on $\mathbb{A}_{e}(-\frac12)$, and equals $-9\sqrt{\theta}\,|y|^{2}$ on $\mathbb{A}_{e}(-1)$.
\end{lemma}
\begin{proof}
At a critical point $x_{0}$ of $\hat{u}$ with Lagrange multiplier $\lambda$, the Hessian of $\hat{u}$ on $T_{x_{0}}\Sigma$ is the restriction of the Hessian of $u-\lambda\langle x,x\rangle $, that is, of $6L(x_{0})-2\lambda\,\operatorname{id}$. Here $2\lambda=3u(x_{0})=3\sqrt{\theta}$ and $L(x_{0})=\sqrt{\theta}L(e)$, so the Hessian is $3\sqrt{\theta}\,(2L(e)-\operatorname{id})=-3\sqrt{\theta}\,\mathcal{K}_{e}$. The values on the Peirce spaces follow from Lemma~\ref{prop:peirce}.
\end{proof}

Note that the Hessian \eqref{eq:hessU} is, up to the factor $-3\sqrt{\theta}$, the operator $\mathcal{K}_{e}$ of \eqref{eq:Ke}: the second variation of $u$ on the sphere and the linearization of the idempotent equation are the same object. Recall that a smooth function on a manifold is \emph{Morse--Bott} along a closed submanifold $M$ of critical points if, at every point of $M$, the kernel of the Hessian equals the tangent space of $M$, and that it is a \emph{Morse--Bott function} if its critical set is a disjoint union of closed submanifolds along each of which it is Morse--Bott \cite{Bott54,BanyagaHurtubise04}.

The critical points with value $0$ are treated in the same way as the idempotents in Section~\ref{sec:euclid}. The role of the defect identity is played by the following consequence of \eqref{eq:V}. If $y\in\mathbb{A}$ and $\delta=y^{2}$, then $y^{3}=y\delta$, $yy^{3}=y(y\delta)$, $y^{2}y^{2}=\delta^{2}$ and $u(y)=\langle y,\delta\rangle $, so that \eqref{eq:V} becomes
\begin{equation}\label{eq:nildefect}
4\,y(y\delta)-3\theta q(y)\,\delta-2\theta\langle y,\delta\rangle \,y=-\delta^{2}.
\end{equation}

\begin{lemma}\label{prop:nilpotent}
Let $x_{0}\in\Sigma$ satisfy $x_{0}^{2}=0$. Then, near $x_{0}$, the set $\operatorname{Nil}(\mathbb{A})=\{x\in\mathbb{A}:x^{2}=0\}$ is a real analytic submanifold of $\mathbb{A}$ with tangent space $\ker L(x_{0})$ at $x_{0}$. Consequently $\operatorname{Nil}_{\Sigma}(\mathbb{A})=\operatorname{Nil}(\mathbb{A})\cap\Sigma$ is a compact real analytic submanifold of $\Sigma$ with tangent space $\ker L(x_{0})\cap x_{0}^{\perp}$ at $x_{0}$, and the Hessian of $\hat{u}$ at $x_{0}$ is the quadratic form $6\langle L(x_{0})y,y\rangle $ on $x_{0}^{\perp}$, whose kernel is this tangent space.
\end{lemma}
\begin{proof}
The argument is that of Lemma~\ref{lem:reduction} and Theorem~\ref{thm:idmsmooth}, with $-\mathcal{K}_{e}$ replaced by $2L(x_{0})$. Put $\mathbb{V}_{0}=\ker L(x_{0})$. Since $L(x_{0})$ is self-adjoint, $\mathbb{V}_{0}^{\perp}$ is its range, $L(x_{0})$ is invertible on $\mathbb{V}_{0}^{\perp}$, and $x_{0}\in\mathbb{V}_{0}$ because $L(x_{0})x_{0}=x_{0}^{2}=0$. For $y=x_{0}+v+w$ with $v\in\mathbb{V}_{0}$, $w\in\mathbb{V}_{0}^{\perp}$ we have
\begin{equation}\label{eq:nilsplit}
y^{2}=2x_{0}w+(v+w)^{2}.
\end{equation}
As in Lemma~\ref{lem:reduction}, the implicit function theorem gives a ball $U\subset\mathbb{V}_{0}$ centered at $0$ and a real analytic map $g:U\to\mathbb{V}_{0}^{\perp}$ with $g(0)=0$, $dg(0)=0$, such that, for $y(v)=x_{0}+v+g(v)$, the element $\delta(v)=y(v)^{2}$ lies in $\mathbb{V}_{0}$, and near $x_{0}$ the $2$-nilpotents are exactly the points $y(v)$ with $\delta(v)=0$. For $v\in U$ define the linear map $\Psi_{v}:\mathbb{V}_{0}\to\mathbb{A}$,
\begin{equation}\label{eq:Psinil}
\Psi_{v}(\eta)=4\,y(y\eta)-3\theta q(y)\,\eta-2\theta\langle y,\eta\rangle \,y,\qquad y=y(v).
\end{equation}
By \eqref{eq:nildefect}, $\Psi_{v}(\delta(v))=-\delta(v)^{2}$. Since $q(x_{0})=1$ and $L(x_{0})\eta=0$ for $\eta\in\mathbb{V}_{0}$, we have $\Psi_{0}(\eta)=-3\theta\eta-2\theta\langle x_{0},\eta\rangle x_{0}$. If $\Psi_{0}(\eta)=0$, pairing with $x_{0}$ gives $-5\theta\langle x_{0},\eta\rangle =0$, hence $\langle x_{0},\eta\rangle =0$ and $\eta=0$; so $\Psi_{0}$ is injective. Now argue as in the proof of Theorem~\ref{thm:idmsmooth}: $\delta(0)=x_{0}^{2}=0$, and if the homogeneous parts $\delta_{j}$ of $\delta$ vanish for $j<k$, then the part of degree $k$ of the identity $\Psi_{v}(\delta(v))=-\delta(v)^{2}$ reads $\Psi_{0}(\delta_{k})=0$, so $\delta_{k}=0$. Hence $\delta$ vanishes to infinite order at $0$, and, being real analytic, $\delta\equiv0$ on $U$. Thus, near $x_{0}$, $\operatorname{Nil}(\mathbb{A})$ is the graph of $g$ over a ball in $\mathbb{V}_{0}$, a real analytic submanifold with tangent space $\mathbb{V}_{0}$ at $x_{0}$.

The set $\operatorname{Nil}(\mathbb{A})$ is a cone, and its tangent space at $x_{0}$ contains the radial direction $x_{0}$. Hence it meets $\Sigma$ transversally, and $\operatorname{Nil}_{\Sigma}(\mathbb{A})$ is a real analytic submanifold of $\Sigma$ with tangent space $\mathbb{V}_{0}\cap x_{0}^{\perp}$; it is compact, being closed in $\Sigma$. Finally, $x_{0}$ is a critical point of $\hat{u}$ with $u(x_{0})=0$, so the Lagrange multiplier vanishes and, as in the proof of Lemma~\ref{prop:hessianU}, the Hessian of $\hat{u}$ at $x_{0}$ is the restriction of $\operatorname{Hess}u(x_{0})=6L(x_{0})$ to $x_{0}^{\perp}$. The subspace $x_{0}^{\perp}$ is $L(x_{0})$-invariant, because $\langle x_{0}y,x_{0}\rangle =\langle y,x_{0}^{2}\rangle =0$, so the kernel of the Hessian is $\ker L(x_{0})\cap x_{0}^{\perp}$.
\end{proof}

\begin{proof}[Proof of Theorem~\ref{mthm:MB}]

(i) By Lemma~\ref{lem:critical}, the critical points of $\hat{u}$ with positive value form the set $M=\sqrt{\theta}\,\operatorname{Idm}(\mathbb{A})$, all with value $\sqrt{\theta}$. The maximum of $\hat{u}$ on the compact sphere is attained at a critical point, and it is positive because $u$ is odd and not identically zero (Lemma~\ref{prop:existence}). Hence $\max_{\Sigma}\hat{u}=\sqrt{\theta}$ and $M$ is the set of maximum points.

(ii) By Lemma~\ref{lem:critical}, the critical set of $\hat{u}$ is the disjoint union of $M$, $-M$ and $\operatorname{Nil}_{\Sigma}(\mathbb{A})$, three closed sets on which $\hat{u}$ takes the values $\sqrt{\theta}$, $-\sqrt{\theta}$ and $0$. By Theorem~\ref{mthm:A}, $M$ is a compact real analytic submanifold of $\Sigma$ with $T_{x_{0}}M=\mathbb{A}_{e}(\frac12)$ at $x_{0}=\sqrt{\theta}e$. By Lemma~\ref{prop:hessianU}, the kernel of $\operatorname{Hess}_{x_{0}}\hat{u}$ is exactly $\mathbb{A}_{e}(\frac12)$, and the Hessian is negative definite on its orthogonal complement $\mathbb{A}_{e}(-1)\oplus\mathbb{A}_{e}(-\frac12)$ in $e^{\perp}$. Since $\hat{u}$ is odd, the same holds at the points of $-M$, with the Hessian changing sign. The set $\operatorname{Nil}_{\Sigma}(\mathbb{A})$ is treated in Lemma~\ref{prop:nilpotent}.

(iii) Let $M_{1},\dots,M_{l}$ be the connected components of $M$; there are finitely many since $M$ is a compact manifold. Let $\Omega=\{x\in \Sigma:\hat{u}(x)>0\}$ and let $\varphi_{t}$ be the flow of the gradient field $\nabla_{\Sigma}\hat{u}$ on the compact manifold $\Sigma$. Along trajectories $\hat{u}$ is nondecreasing, so $\Omega$ is forward invariant.

By (ii) and the Morse--Bott lemma \cite{Bott54,BanyagaHurtubise04}, near $M_{i}$ there are coordinates in which $\hat{u}=\sqrt{\theta}-|\nu|^{2}$, where $\nu$ is a coordinate in the normal directions. Since $M$ is compact, these coordinates can be chosen on a tubular neighborhood $T$ of $M$, and since $\hat{u}<\sqrt{\theta}$ on the compact set $\Sigma\setminus T$, the set $\{\hat{u}>\sqrt{\theta}-\varepsilon\}$ lies in $T$ for small $\varepsilon>0$. Hence, for $\varepsilon>0$ small enough, the superlevel set $\{\hat{u}>\sqrt{\theta}-\varepsilon\}$ is the disjoint union of connected open sets $O_{1},\dots,O_{l}\subset\Omega$ with $M_{i}\subset O_{i}$. Since $\hat{u}$ is nondecreasing along trajectories and trajectories are connected, a trajectory that enters $O_{i}$ stays in $O_{i}$. Put
\begin{equation*}
B_{i}=\{x\in\Omega:\ \varphi_{t}(x)\in O_{i}\text{ for some }t\geq0\}.
\end{equation*}
Each $B_{i}$ is open, being the union of the open sets $\varphi_{t}^{-1}(O_{i})$, and $B_{i}\cap B_{j}=\emptyset$ for $i\neq j$.

Let $x\in\Omega$. Its $\omega$-limit set is nonempty, consists of critical points of $\hat{u}$, and lies in $\{\hat{u}\geq \hat{u}(x)\}$ with $\hat{u}(x)>0$; by Lemma~\ref{lem:critical} it lies in $M$. Since $\hat{u}$ tends to $\sqrt{\theta}$ along the trajectory, the trajectory enters $\{\hat{u}>\sqrt{\theta}-\varepsilon\}$, hence some $O_{i}$. So $x\in B_{i}$, and $\Omega$ is the disjoint union of the open sets $B_{1},\dots,B_{l}$. Each $B_{i}$ is connected, because every point of $B_{i}$ is joined by its trajectory to the connected set $O_{i}$. Hence the components of $\Omega$ are exactly $B_{1},\dots,B_{l}$, and $B_{i}\cap M=M_{i}$, since $M\cap O_{i}=M_{i}$ and $M$ consists of fixed points of the flow.
\end{proof}

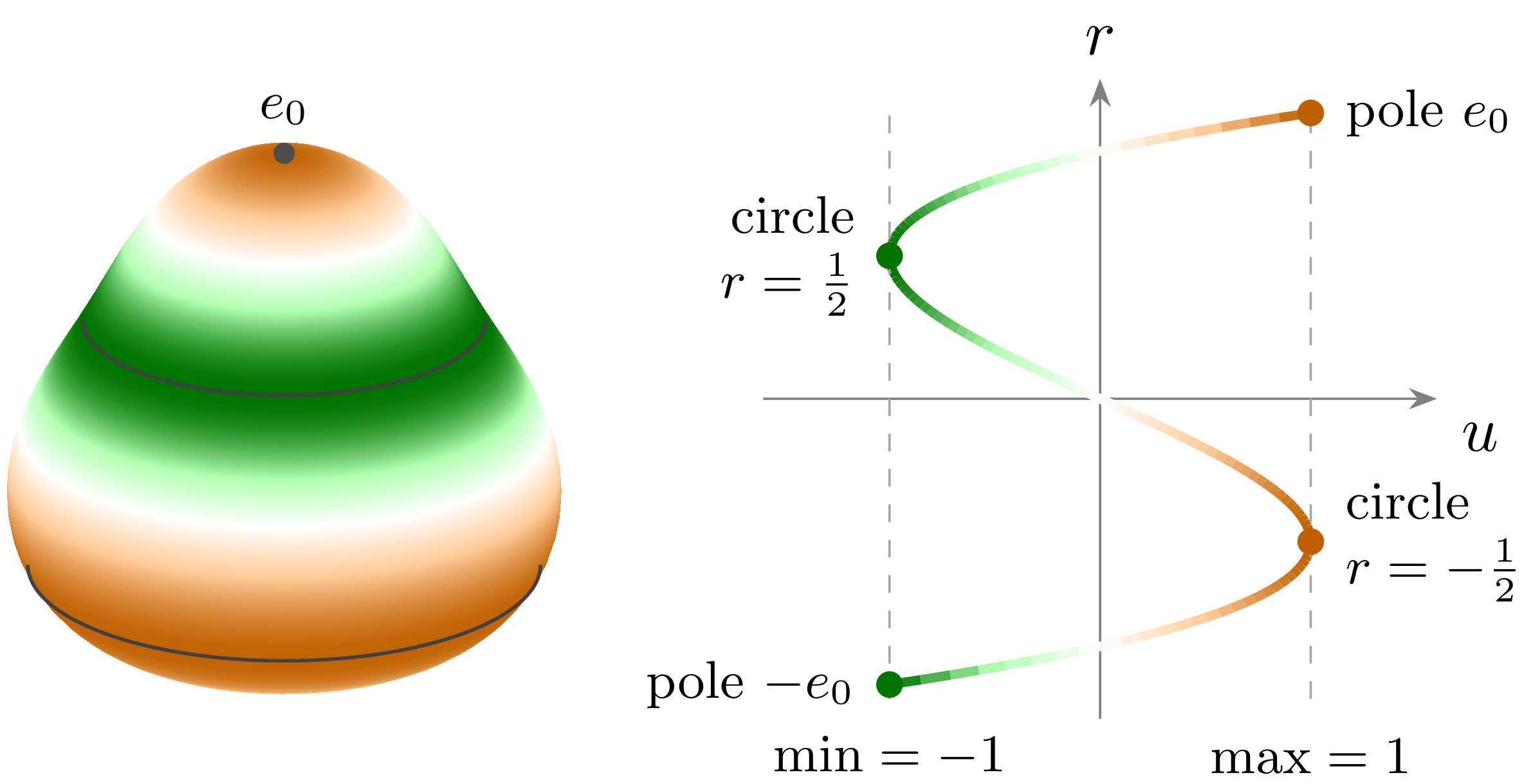
\begin{figure}[ht]
\centering
\resizebox{0.92\textwidth}{!}{%
\begin{tikzpicture}
\begin{axis}[name=sph,
  width=7.2cm,height=7.2cm,axis equal image,hide axis,
  view={30}{22},
  colormap={orgr}{color(0cm)=(green!45!black) color(1cm)=(green!30) color(2cm)=(white) color(3cm)=(orange!40) color(4cm)=(orange!75!black)},
  point meta min=-1,point meta max=1,
  z buffer=sort]
\addplot3[surf,shader=interp,variable=\t,variable y=\s,domain=0:180,y domain=0:360,samples=46,samples y=72,
  point meta={4*cos(\t)^3-3*cos(\t)}]
  ({(1+0.12*(4*cos(\t)^3-3*cos(\t)))*sin(\t)*cos(\s)},
   {(1+0.12*(4*cos(\t)^3-3*cos(\t)))*sin(\t)*sin(\s)},
   {(1+0.12*(4*cos(\t)^3-3*cos(\t)))*cos(\t)});
\draw[black!75,line width=0.6pt] plot[variable=\s,domain=-150:30,samples=60,smooth] (axis cs:{0.88*sin(60)*cos(\s)},{0.88*sin(60)*sin(\s)},{0.88*cos(60)});
\draw[black!75,line width=0.6pt] plot[variable=\s,domain=-150:30,samples=60,smooth] (axis cs:{1.12*sin(120)*cos(\s)},{1.12*sin(120)*sin(\s)},{1.12*cos(120)});
\addplot3[only marks,mark=*,mark size=1.6pt,black!70,samples y=1] coordinates {(0,0,1.12)};
\node[font=\scriptsize,anchor=south] at (axis cs:0,0,1.14) {$e_{0}$};
\end{axis}
\begin{axis}[at={(sph.east)},anchor=west,xshift=0.6cm,
  width=5.6cm,height=5.4cm,
  xmin=-1.6,xmax=1.6,ymin=-1.12,ymax=1.12,
  axis lines=middle,axis line style={black!50,-{Stealth}},
  xtick=\empty,ytick=\empty,
  xlabel={$u$},ylabel={$r$},
  xlabel style={anchor=north west,font=\small},ylabel style={anchor=south,font=\small},
  colormap={orgr}{color(0cm)=(green!45!black) color(1cm)=(green!30) color(2cm)=(white) color(3cm)=(orange!40) color(4cm)=(orange!75!black)},
  point meta min=-1,point meta max=1,clip=false]
\addplot[mesh,line width=1.6pt,domain=-1:1,samples=120,point meta=x]
  ({4*x^3-3*x},{x});
\draw[black!35,dashed] (axis cs:1,-1.05)--(axis cs:1,1.05);
\draw[black!35,dashed] (axis cs:-1,-1.05)--(axis cs:-1,1.05);
\node[circle,fill=orange!75!black,inner sep=1.6pt] at (axis cs:1,1) {};
\node[circle,fill=orange!75!black,inner sep=1.6pt] at (axis cs:1,-0.5) {};
\node[circle,fill=green!45!black,inner sep=1.6pt] at (axis cs:-1,0.5) {};
\node[circle,fill=green!45!black,inner sep=1.6pt] at (axis cs:-1,-1) {};
\node[font=\scriptsize,anchor=west,align=left] at (axis cs:1.05,1) {pole $e_{0}$};
\node[font=\scriptsize,anchor=west,align=left] at (axis cs:1.05,-0.5) {circle\\ $r=-\frac12$};
\node[font=\scriptsize,anchor=east,align=right] at (axis cs:-1.05,0.5) {circle\\ $r=\frac12$};
\node[font=\scriptsize,anchor=east] at (axis cs:-1.05,-1) {pole $-e_{0}$};
\node[font=\scriptsize] at (axis cs:1,-1.25) {$\max=1$};
\node[font=\scriptsize] at (axis cs:-1,-1.25) {$\min=-1$};
\end{axis}
\end{tikzpicture}%
}
\caption{Example~\ref{ex:rpc} with $m=2$, so that $\mathbb{A}=\mathbb{R}^{3}$ and $\theta=1$. Left: the radial graph $\omega\mapsto(1+0.12\,u(\omega))\,\omega$ over the unit sphere $\Sigma$, colored by $u$ (orange for $u>0$, green for $u<0$). Right: on $\Sigma$, $u$ is the Chebyshev polynomial $T_{3}(r)=4r^{3}-3r$ of the height $r$. The maximum set $\sqrt{\theta}\,\operatorname{Idm}(\mathbb{A})$ consists of the pole $e_{0}$ (with $n_{1/2}=0$) and the circle $r=-\frac12$ (with $n_{1/2}=1$); the minimum set consists of their antipodes. These are the only critical points, since $\mathbb{A}$ has no $2$-nilpotents, and $u|_{\Sigma}$ is Morse--Bott along each of them.}
\label{fig:chebyshev}
\end{figure}

\begin{remark}
By Theorem~\ref{mthm:MB}(ii), each $M_{i}$ is a normally hyperbolic attracting manifold of fixed points of the gradient flow. By the theory of normally hyperbolic invariant manifolds \cite{HirschPughShub77}, its basin $B_{i}$ is fibered over $M_{i}$ by the stable manifolds of its points, which are open discs of dimension $n_{-1}+n_{-1/2}$; in particular $B_{i}$ deformation retracts onto $M_{i}$. In Example~\ref{ex:rpc}, with $m\geq2$ and $z=(r,x)$, one has
\begin{equation*}
u(z)=r^{3}-3r|x|^{2},
\end{equation*}
and on $\Sigma$ the region $u>0$ consists of the cap $r>\frac{\sqrt3}{2}$, a disc around $e_{0}$, and the band $-\frac{\sqrt3}{2}<r<0$, which retracts onto the sphere of idempotents $(-\frac12,v)$. For $m=2$ this is shown in Figure~\ref{fig:chebyshev}: on $\Sigma=S^{2}$ the cubic form depends only on the height $r$ and equals the Chebyshev polynomial $T_{3}(r)=4r^{3}-3r$, so that its critical sets are two points and two circles, in accordance with Theorem~\ref{mthm:MB}; compare the role of equioscillation in \cite{Tkachev26}. In Example~\ref{ex:rac}, $u=3\det$, and the region $\det x>0$ of the unit sphere has two components: the matrices with three positive eigenvalues, retracting onto $I/\sqrt3$, and those with one positive and two negative eigenvalues, retracting onto the projective plane $\{(2ww^{T}-I)/\sqrt3\}$.
\end{remark}

\begin{remark}
Statement (i) of Theorem~\ref{mthm:MB} is a rigidity property: the cubic form $u$ has no positive local maxima, saddle points or minima on $\Sigma$ other than its global maximum set. Critical points with value $0$ are the normalized $2$-nilpotents ($x^{2}=0$), and the negative critical points are $-M$. Since $\hat{u}$ is real analytic, the {\L}ojasiewicz gradient inequality \cite{Lojasiewicz63} implies that every trajectory of $\nabla_{\Sigma}\hat{u}$ in $\Omega$ converges to a single point of $M$. In terms of the flow \eqref{eq:flow}, every trajectory starting in the cone over $\{u>0\}$ is therefore asymptotic in direction to an idempotent: its radial projection to $\Sigma$ converges to a point of $\sqrt{\theta}\,\operatorname{Idm}(\mathbb{A})$.
\end{remark}

\section{The flow \texorpdfstring{$\dot x=x^{2}$}{x\'=x\^{}2} and two first integrals}\label{sec:flow}

Why should the identity \eqref{eq:hsiang} produce conserved quantities? The answer is a classical computation from the theory of the $\infty$-Laplacian \cite{Aronsson67,Crandall08}. Along a trajectory of a gradient flow $\dot x=\nabla v$ one has
\begin{equation}\label{eq:aronsson}
\frac{d}{dt}\,|\nabla v|^{2}=2\langle (\operatorname{Hess}v)\nabla v,\nabla v\rangle =2\Delta_{\infty}v ,
\end{equation}
which is the basic reason why $\Delta_{\infty}$ governs the behavior of a function along its gradient curves. For $v=\frac13u$ the flow is \eqref{eq:flow}, and \eqref{eq:aronsson} together with \eqref{eq:radialinf} gives
\begin{equation}\label{eq:speed}
\frac{d}{dt}\,|\dot x|^{2}=4\theta\,|x|^{2}u(x)=\theta\,\frac{d}{dt}\,|x|^{4}.
\end{equation}
Thus \eqref{eq:radialinf} says precisely that the squared speed $|\dot x|^{2}$ and the quartic quantity $\theta|x|^{4}$ change at the same rate, that is, that the energy $\mathcal{E}=|\dot x|^{2}-\theta|x|^{4}$ is conserved. Radial $\infty$-Hsiang algebras are therefore exactly the quadratic gradient systems with this quartic energy, which is a further natural source of these algebras besides minimal cones and nonlinear elliptic equations. The second conserved quantity $\mathcal{F}$ is less obvious; it comes from a hidden linear relation between two invariants of degree $7$, see \eqref{eq:AB1} and \eqref{eq:AB2}. The two integrals are also the key to the rigidity statements of Subsection~\ref{ssec:1.2}: their gradients are, up to a factor, the defects $\mathcal{R}$ and $\mathcal{S}$, so the vector identities $\mathcal{R}=0$ and $\mathcal{S}=0$ become the identical vanishing of global scalar quantities. Below we carry out the computations algebraically, replacing the time derivative by a derivation.

For a polynomial $f$ on $\mathbb{A}$ define
\begin{equation}\label{eq:Lie}
(\partial f)(x)=df_{x}(x^{2})=\langle \nabla f(x),x^{2}\rangle ,
\end{equation}
the derivative of $f$ along solutions of \eqref{eq:flow}. We call $f$ a \emph{first integral} if $\partial f\equiv0$, that is, if $f$ is constant along all solutions. The operator $\partial$ is a derivation, and, applying it componentwise to vector-valued polynomials, we get
\begin{equation}\label{eq:flowrules}
\partial(x^{2})=2x^{3},\quad \partial(x^{3})=2xx^{3}+x^{2}x^{2},\quad \partial q=2u,\quad \partial u=3\langle x^{2},x^{2}\rangle .
\end{equation}
For example, $\partial u=\langle x^{2},x^{2}\rangle +2\langle x,x^{3}\rangle $, and $\langle x,x^{3}\rangle =\langle x^{2},x^{2}\rangle $ by invariance. Recall the energy and the sextic integral from \eqref{eq:EFintro}:
\begin{equation}\label{eq:EF}
\begin{aligned}
\mathcal{E}(x)&=\langle x^{2},x^{2}\rangle -\theta q(x)^{2},\\
\mathcal{F}(x)&=\langle x^{2},x^{2}x^{2}\rangle -\theta u(x)^{2}=\langle \mathcal{S}(x),x^{2}\rangle .
\end{aligned}
\end{equation}
Along a solution of \eqref{eq:flow}, $\mathcal{E}=|\dot x|^{2}-\theta|x|^{4}$, which explains the name.

\begin{theorem}\label{thm:integrals}
Let $\mathbb{A}$ be a commutative algebra with a metric and $\theta\in\mathbb{R}^{\times}$.
\begin{itemize}
\item[{(i)}] $\nabla \mathcal{E}=4\mathcal{R}$ and $\partial \mathcal{E}=4\mathcal{P}$, with $\mathcal{P}$ as in \eqref{eq:P}. Hence $\mathbb{A}$ is radial $\infty$-Hsiang with constant $\theta$ if and only if $\mathcal{E}$ is a first integral.
\item[{(ii)}] $\nabla \mathcal{F}=6L(x)\mathcal{S}(x)$, and if $\mathbb{A}$ is radial $\infty$-Hsiang, then $\mathcal{F}$ is a first integral.
\item[{(iii)}] $\mathcal{E}\equiv0$ if and only if $\mathbb{A}$ is radial pseudocomposition.
\item[{(iv)}] If $\mathbb{A}$ is radial $\infty$-Hsiang and $\operatorname{Idm}(\mathbb{A})\neq\emptyset$, then $\mathcal{F}\equiv0$ if and only if $\mathbb{A}$ is radial admissible cubic.
\end{itemize}
\end{theorem}
\begin{proof}
(i) Since $d\langle x^{2},x^{2}\rangle _{x}(y)=4\langle x^{2},xy\rangle =\langle 4x^{3},y\rangle $ and $d(q^{2})_{x}(y)=4q\langle x,y\rangle $,
\begin{equation}\label{eq:gradE}
\begin{aligned}
\nabla \mathcal{E}&=4\bigl(x^{3}-\theta q\,x\bigr)=4\mathcal{R},\\
\partial \mathcal{E}&=\langle 4\mathcal{R},x^{2}\rangle =4\mathcal{P},
\end{aligned}
\end{equation}
the last equality by \eqref{eq:shadow}.

(ii) The cubic form $w\mapsto\langle w,w^{2}\rangle $ has differential $3\langle w^{2},\cdot\rangle $. With $w=x^{2}$ and $dw_{x}(y)=2xy$ we get
\begin{equation}\label{eq:dF}
\begin{aligned}
d\langle x^{2},x^{2}x^{2}\rangle _{x}(y)&=6\langle x^{2}x^{2},xy\rangle =6\langle x(x^{2}x^{2}),y\rangle ,\\
d(u^{2})_{x}(y)&=6u\langle x^{2},y\rangle .
\end{aligned}
\end{equation}
Hence $\nabla \mathcal{F}=6\bigl(x(x^{2}x^{2})-\theta u\,x^{2}\bigr)=6L(x)\mathcal{S}(x)$.

Assume now \eqref{eq:hsiang}, and put
\begin{equation}\label{eq:AB}
\begin{aligned}
A&=\langle x^{3},xx^{3}\rangle ,\\
B&=\langle x^{3},x^{2}x^{2}\rangle .
\end{aligned}
\end{equation}
Since $\partial \mathcal{F}=\langle \nabla \mathcal{F},x^{2}\rangle =6\langle x^{2}x^{2},x^{3}\rangle -6\theta u\langle x^{2},x^{2}\rangle $, we must show that $B=\theta u\langle x^{2},x^{2}\rangle $. We obtain two linear relations between $A$ and $B$. Pairing \eqref{eq:V} with $x^{3}$, and using $\langle x,x^{3}\rangle =\langle x^{2},x^{2}\rangle $ and $\langle x^{2},x^{3}\rangle =\theta qu$, gives
\begin{equation}\label{eq:AB1}
4A+B=2\theta u\langle x^{2},x^{2}\rangle +3\theta^{2}q^{2}u .
\end{equation}
Pairing \eqref{eq:V} with $x^{2}$, and using $\langle xx^{3},x^{2}\rangle =\langle x^{3},x^{3}\rangle $, gives the identity
\begin{equation}\label{eq:H}
\mathcal{Z}(x):=4\langle x^{3},x^{3}\rangle +\langle x^{2},x^{2}x^{2}\rangle -2\theta u^{2}-3\theta q\langle x^{2},x^{2}\rangle \equiv0,
\end{equation}
so that $\partial \mathcal{Z}\equiv0$. By \eqref{eq:flowrules} and invariance,
\begin{equation}\label{eq:Drules}
\begin{aligned}
\partial\langle x^{3},x^{3}\rangle &=4A+2B,\\
\partial\langle x^{2},x^{2}x^{2}\rangle &=6B,\\
\partial(u^{2})&=6u\langle x^{2},x^{2}\rangle ,\\
\partial\langle x^{2},x^{2}\rangle &=4\theta qu,\\
\partial\bigl(q\langle x^{2},x^{2}\rangle \bigr)&=2u\langle x^{2},x^{2}\rangle +4\theta q^{2}u .
\end{aligned}
\end{equation}
Hence
\begin{equation}\label{eq:AB2}
0=\partial \mathcal{Z}=16A+14B-18\theta u\langle x^{2},x^{2}\rangle -12\theta^{2}q^{2}u .
\end{equation}
Subtracting four times \eqref{eq:AB1} from \eqref{eq:AB2} gives $10B=10\theta u\langle x^{2},x^{2}\rangle $, so $B=\theta u\langle x^{2},x^{2}\rangle $, and therefore $\partial \mathcal{F}=0$.

(iii) If $\mathcal{E}\equiv0$, then $\mathcal{R}=\frac14\nabla \mathcal{E}\equiv0$. Conversely, if $\mathcal{R}\equiv0$, then $\mathcal{E}(x)=\langle x^{3},x\rangle -\theta q^{2}=\langle \mathcal{R}(x),x\rangle =0$.

(iv) If $\mathcal{S}\equiv0$, then $\mathcal{F}=\langle \mathcal{S},x^{2}\rangle \equiv0$. Conversely, if $\mathcal{F}\equiv0$, then $L(x)\mathcal{S}(x)\equiv0$, and multiplying by the adjugate of $L(x)$ gives $\det L(x)\,\mathcal{S}(x)\equiv0$. For $e\in\operatorname{Idm}(\mathbb{A})$, $L(e)$ has no zero eigenvalue by Lemma~\ref{prop:peirce}, so $\det L(x)$ is a nonzero polynomial and $\mathcal{S}$ vanishes on the nonempty open set where $\det L(x)\neq0$. Hence $\mathcal{S}\equiv0$.
\end{proof}

Theorem~\ref{thm:integrals}(i),(ii) gives Theorem~\ref{mthm:B}(i) and the gradient formulas in Theorem~\ref{mthm:B}(ii), and Theorem~\ref{thm:integrals}(iii),(iv) gives the rest of Theorem~\ref{mthm:B}(ii).

\begin{remark}\label{rem:EKP}
The sextic integral $\mathcal{F}$ has a transparent meaning in terms of the gradient map of $u$. By \eqref{eq:dictionary}, $\Gamma(x)=x^{2}=\frac13\nabla u(x)$, and $\mathcal{F}(x)=u(\Gamma(x))-\theta\,u(x)^{2}$. Thus $\mathcal{F}\equiv0$ means that the cubic form $u$ is multiplicative on squares,
\begin{equation}\label{eq:umult}
u(x^{2})=\theta\,u(x)^{2},
\end{equation}
just as $\mathcal{E}\equiv0$ means that $|x^{2}|^{2}=\theta|x|^{4}$. This is the radial counterpart of a basic property of admissible cubic algebras: if $x^{2}x^{2}=N(x)x$, then applying this identity to $x^{2}$ and using $(x^{2})^{2}=N(x)x$ gives $N(x^{2})\,x^{2}=N(x)^{2}x^{2}$, and, if $N\not\equiv0$, a density argument yields, writing $x^{\#}=x^{2}$ as in the adjoint identity $(x^{\#})^{\#}=N(x)x$ of cubic Jordan algebras,
\begin{equation}\label{eq:Nmult}
N(x^{\#})=N(x)^{2}\qquad\text{for all }x\in\mathbb{A};
\end{equation}
see \cite[Section~2, Eq.~(2)]{ElduqueOkubo00}. For the determinant in Example~\ref{ex:rac}, \eqref{eq:Nmult} is the classical identity $\det(x^{\#})=\det(x)^{2}$. In a radial admissible cubic algebra the identity $x^{2}x^{2}=\theta u(x)\,x$ says moreover that
\begin{equation}\label{eq:cremona}
\Gamma(\Gamma(x))=\theta\,u(x)\,x .
\end{equation}
Hence the gradient map of $u$ induces a birational involution of the projective space $P(\mathbb{A})$, a quadro-quadric Cremona transformation. Identifying $\mathbb{A}$ with its dual by the metric, the multiplicative Legendre transform $u_{*}$, defined by $u_{*}(\nabla u(x)/u(x))=1/u(x)$, is then again a cubic form, namely $u_{*}=u/(27\theta)$. Cubic forms with this property were studied by Etingof, Kazhdan and Polishchuk \cite{EtingofKazhdanPolishchuk02} in connection with the question of when the Fourier transform of an elementary function is elementary: for the phase $f(x)/t$ with a cubic form $f$, the corresponding semiclassical condition holds precisely for this class. Over $\mathbb{C}$ they showed that such cubic forms are given by a short explicit list, essentially the generic norms of semisimple Jordan algebras of degree~$3$, such as the determinant in Example~\ref{ex:rac}; the proof uses Zak's classification of Severi varieties. In our notation, the singular locus of the cubic hypersurface $\{u=0\}$ is the cone $\operatorname{Nil}(\mathbb{A})=\{x:x^{2}=0\}$ of $2$-nilpotents, whose projectivization plays the role of the Severi variety in that argument; in Example~\ref{ex:rac} it is the cone over the Veronese surface of rank-one symmetric matrices, and its smoothness is consistent with Lemma~\ref{prop:nilpotent}. Theorem~\ref{mthm:B} adds a dynamical characterization of this class: it consists of the radial $\infty$-Hsiang algebras with an idempotent for which the first integral $\mathcal{F}$ of the flow \eqref{eq:flow} vanishes identically.
\end{remark}

\begin{remark}
The divergence of the vector field $x\mapsto x^{2}$ is $2\operatorname{tr} L(x)$. So $\mathbb{A}$ is \emph{exact}, meaning $\operatorname{tr} L(x)\equiv0$, if and only if the flow \eqref{eq:flow} preserves volume. In the forthcoming paper mentioned in Subsection~\ref{ssec:1.2}, we prove that a radial $\infty$-Hsiang algebra all of whose idempotents, including those of its complexification, have maximal Peirce complexity is exact; for such algebras the flow \eqref{eq:flow} therefore preserves volume. This result is not used in the present paper.
\end{remark}

\section{Idempotents as Darboux points}\label{sec:darboux}

For $e\in\operatorname{Idm}(\mathbb{A})$ the curve $x(t)=-t^{-1}e$ solves \eqref{eq:flow}, because
\begin{equation*}
\frac{d}{dt}\bigl(-t^{-1}e\bigr)=t^{-2}e=\bigl(-t^{-1}e\bigr)^{2}.
\end{equation*}
So $c=-e$ is a Darboux point in the sense of Subsection~\ref{ssec:1.4} and of \cite{Yoshida83,Furta96,Kozlov96,MaciejewskiPrzybylska05}. Conversely, every Darboux point $c$ ($c^{2}=-c\neq0$) is of this form, with $e=-c$. The Kovalevskaya matrix (see \cite{Goriely01}) is
\begin{equation}\label{eq:Kmatrix}
\mathcal{K}(-e)=\operatorname{id}+2L(-e)=\operatorname{id}-2L(e)=\mathcal{K}_{e},
\end{equation}
with eigenvalues \eqref{eq:kovalevskaya}. Explicitly, by \eqref{eq:variational}, for $\xi_{0}\in\mathbb{A}_{e}(\mu)$ the function $\xi(t)=t^{-2\mu}\xi_{0}=t^{\rho-1}\xi_{0}$ solves the variational equation along $x(t)=-t^{-1}e$. Thus perturbations in $\mathbb{A}_{e}(-1)$ and $\mathbb{A}_{e}(-\frac12)$ grow like $t^{2}$ and $t$, while perturbations in $\mathbb{A}_{e}(\frac12)$ decay like $t^{-1}$, at the same rate as the solution itself. The exponent $0$ belongs exactly to the tangent space $\mathbb{A}_{e}(\frac12)$ of the idempotent manifold (Theorem~\ref{thm:idmsmooth}): the Darboux points form the set $-\operatorname{Idm}(\mathbb{A})$, and the zero exponents account for the directions along this family. All Kovalevskaya exponents are integers, namely $-1,0,2,3$.

\begin{theorem}\label{thm:hessians}
Let $\mathbb{A}$ be radial $\infty$-Hsiang and $e\in\operatorname{Idm}(\mathbb{A})$. Then the first integrals $\mathcal{E}$ and $\mathcal{F}$ defined in \eqref{eq:EF} satisfy $\mathcal{E}(e)=\mathcal{F}(e)=0$ and $\nabla \mathcal{E}(e)=\nabla \mathcal{F}(e)=0$. For $y\in\mathbb{A}_{e}(\mu)$ and $z\in\mathbb{A}_{e}(\mu')$ with $\mu,\mu'\neq1$,
\begin{equation}\label{eq:hessians}
\begin{aligned}
d^{2}\mathcal{E}_{e}(y,z)&=4(2\mu-1)(\mu+1)\langle y,z\rangle ,\\
d^{2}\mathcal{F}_{e}(y,z)&=6\mu(2\mu-1)(2\mu+1)\langle y,z\rangle .
\end{aligned}
\end{equation}
Consequently
\begin{equation}\label{eq:ranks}
\begin{aligned}
\operatorname{rank}\bigl(d^{2}\mathcal{E}_{e}|_{e^{\perp}}\bigr)&=n_{-1/2}(e),\\
\operatorname{rank}\bigl(d^{2}\mathcal{F}_{e}|_{e^{\perp}}\bigr)&=n_{-1}(e).
\end{aligned}
\end{equation}
\end{theorem}
\begin{proof}
By Lemma~\ref{lem:defects}, $\nabla \mathcal{E}(e)=4\mathcal{R}(e)=0$ and $\nabla \mathcal{F}(e)=6L(e)\mathcal{S}(e)=0$; also $\mathcal{E}(e)=\mathcal{F}(e)=\langle e,e\rangle -\theta\langle e,e\rangle ^{2}=0$. Differentiating $\nabla \mathcal{E}=4\mathcal{R}$ and $\nabla \mathcal{F}=6L(x)\mathcal{S}(x)$ at $e$, and using $\mathcal{S}(e)=0$, we obtain
\begin{equation}\label{eq:d2EF}
\begin{aligned}
d^{2}\mathcal{E}_{e}(y,z)&=4\langle d\mathcal{R}_{e}(y),z\rangle ,\\
d^{2}\mathcal{F}_{e}(y,z)&=6\langle L(e)\,d\mathcal{S}_{e}(y),z\rangle .
\end{aligned}
\end{equation}
Now \eqref{eq:hessians} follows from \eqref{eq:dRdS}, and both expressions vanish for $\mu\neq\mu'$ by orthogonality of the Peirce spaces. On the orthogonal sum $e^{\perp}=\mathbb{A}_{e}(-1)\oplus\mathbb{A}_{e}(-\frac12)\oplus\mathbb{A}_{e}(\frac12)$, the form $d^{2}\mathcal{E}_{e}$ is $-4\langle \cdot,\cdot\rangle $ on $\mathbb{A}_{e}(-\frac12)$ and zero on the other summands, and $d^{2}\mathcal{F}_{e}$ is $-18\langle \cdot,\cdot\rangle $ on $\mathbb{A}_{e}(-1)$ and zero on the other summands. Since the metric is nondegenerate on each summand, the ranks are as stated.
\end{proof}

This proves Theorem~\ref{mthm:B}(iii). Theorems~\ref{thm:idmsmooth} and \ref{thm:hessians} give a complete geometric reading of the Peirce numbers: $n_{1/2}(e)$ is the dimension of the idempotent manifold at $e$, while $n_{-1/2}(e)$ and $n_{-1}(e)$ measure the nondegeneracy at $e$ of the two global first integrals.

We close this section with the observation mentioned in Subsection~\ref{ssec:1.4}. Along with the first-order flow \eqref{eq:flow}, the cubic form $u$ defines the Newtonian system
\begin{equation}\label{eq:newton}
\ddot z=z^{2}.
\end{equation}
By \eqref{eq:gradqu}, it can be written as
\begin{equation}\label{eq:newtonV}
\ddot z=-\nabla V(z),\qquad V(z)=-\tfrac13u(z),
\end{equation}
with Hamiltonian $H(z,p)=\frac12|p|^{2}+V(z)$ and a homogeneous potential of degree $3$. For such systems, the strongest known necessary conditions for integrability are given by the Morales-Ruiz--Ramis theory \cite{MoralesRuiz99}, and they are expressed through the eigenvalues of the Hessian of $V$ at the Darboux points of the potential, that is, the points $c\neq0$ with $\nabla V(c)=c$ \cite{MoralesRamis01,MaciejewskiPrzybylska04,MaciejewskiPrzybylska05}.

\begin{lemma}\label{prop:newton}
Let $\mathbb{A}$ be a Euclidean radial $\infty$-Hsiang algebra with nonzero multiplication and $\dim\mathbb{A}\geq2$. Then the Hamiltonian system \eqref{eq:newton} is not integrable in the Liouville sense with first integrals that are meromorphic on the complexified phase space.
\end{lemma}
\begin{proof}
We use orthonormal coordinates, so that the kinetic energy is $\frac12\sum p_{i}^{2}$. For $e\in\operatorname{Idm}(\mathbb{A})$ and $c=-e$ we have $\nabla V(c)=-c^{2}=-e=c$, so $c$ is a Darboux point of $V$, and
\begin{equation*}
\operatorname{Hess}V(c)=-2L(c)=2L(e).
\end{equation*}
Its eigenvalues are $2$ on $\mathbb{R}e$ (this is the trivial eigenvalue $k-1$, $k=\deg V=3$) and $\lambda=2\mu$ on $\mathbb{A}_{e}(\mu)$. By the Morales-Ruiz--Ramis theorem \cite{MoralesRamis01}, in the form of \cite[Theorem~1.1]{MaciejewskiPrzybylska05}, if the system is integrable, then for $k=3$ every nontrivial eigenvalue $\lambda$ belongs to one of the sets
\begin{equation*}
\Bigl\{\tfrac{p(3p-1)}{2}\Bigr\},\qquad \Bigl\{\tfrac13+\tfrac{3p(p+1)}{2}\Bigr\},\qquad \Bigl\{-\tfrac1{24}+\alpha\,(\beta+\gamma p)^{2}\Bigr\},\qquad p\in\mathbb{Z},
\end{equation*}
where the last type stands for four families with explicit rational constants $\alpha>0$, $\beta$, $\gamma$. All these numbers are nonnegative, so the eigenvalues $\lambda=-2$ ($\mu=-1$) and $\lambda=-1$ ($\mu=-\frac12$) are forbidden. It therefore suffices to find $e\in\operatorname{Idm}(\mathbb{A})$ with $n_{-1}(e)+n_{-1/2}(e)\geq1$.

Suppose, to the contrary, that $e^{\perp}=\mathbb{A}_{e}(\frac12)$ for every $e\in\operatorname{Idm}(\mathbb{A})$, and fix one such $e$ (it exists by Lemma~\ref{prop:existence}). By Theorem~\ref{mthm:A}, the component of $\operatorname{Idm}(\mathbb{A})$ through $e$ is then a compact submanifold of the sphere \eqref{eq:sphere} of the same dimension $\dim\mathbb{A}-1$. It is therefore open and closed in the sphere, which is connected since $\dim\mathbb{A}\geq2$; so it is the whole sphere. In particular $-e\in\operatorname{Idm}(\mathbb{A})$, which is absurd because $(-e)^{2}=e\neq-e$.
\end{proof}

\begin{remark}
The argument only uses the existence of one idempotent with $n_{-1}+n_{-1/2}\geq1$. Note that the Newtonian system \eqref{eq:newton} and the first-order system \eqref{eq:flow} are different systems; for the latter, all Kovalevskaya exponents \eqref{eq:kovalevskaya} are integers, so Kovalevskaya's test gives no obstruction, and Theorem~\ref{mthm:B} provides two polynomial first integrals.
\end{remark}

\section{Resonances and pseudocomposition algebras}\label{sec:resonance}

In this section we turn the observations of Section~\ref{sec:darboux} into a proof of Theorem~\ref{mthm:C}. By Theorem~\ref{thm:hessians}, at an idempotent $e$ the Hessian of the energy $\mathcal{E}$ of \eqref{eq:EF} sees only $\mathbb{A}_{e}(-\frac12)$. If $n_{-1/2}(e)=0$, this Hessian vanishes, but this alone does not imply $\mathcal{E}\equiv0$, since higher-order terms of $\mathcal{E}$ at $e$ might survive. These terms are controlled by the resonance lemma below: the expansion of any homogeneous first integral at a Darboux point starts with resonant terms. A count of resonances then shows that $\mathcal{E}$ has no such terms when $n_{-1/2}(e)=0$.

The following lemma is an algebraic version of a principle going back to Yoshida and Kozlov \cite{Yoshida83,Kozlov96,Furta96}: near a Darboux point, a homogeneous first integral can only contain terms compatible with the Kovalevskaya exponents. The idea is simple. Write a solution near the special one as $x(t)=-t^{-1}\bigl(e+z(t)\bigr)$ with $z(t)\perp e$. To first order, by \eqref{eq:variational}, the component of $z$ in an eigendirection with exponent $\rho_{i}$ behaves like $t^{\rho_{i}}$ (the corresponding component of the perturbation $-t^{-1}z$ behaves like $t^{\rho_{i}-1}$, as for the solutions of the variational equation), so a monomial $z^{\alpha}$ behaves like $t^{\sum\alpha_{i}\rho_{i}}$. On the other hand, if $I$ is a homogeneous first integral of degree $d$, then $I(e+z(t))=(-t)^{d}I(x(t))$ behaves like $t^{d}$. Only monomials with matching growth can survive in the lowest-order part of $I(e+z)$.

\begin{lemma}[Resonance lemma]\label{lem:resonance}
Let $\mathbb{A}$ be radial $\infty$-Hsiang, $e\in\operatorname{Idm}(\mathbb{A})$, $W=e^{\perp}$, and let $\mathcal{K}_{e}$ be as in \eqref{eq:Ke}; note that $W$ is $\mathcal{K}_{e}$-invariant. Let $I$ be a homogeneous polynomial of degree $d\geq1$ with $\partial I\equiv0$, and put $\psi(z)=I(e+z)$ for $z\in W$. If $\psi\not\equiv0$, let $\psi_{m}$ be its nonzero homogeneous component of lowest degree $m$. Then
\begin{equation}\label{eq:resonance}
d\psi_{m}(z)(\mathcal{K}_{e}z)=d\cdot \psi_{m}(z)\qquad\text{for all }z\in W .
\end{equation}
Equivalently, in Peirce coordinates $z_{1},\dots,z_{N}$ at $e$, that is, in coordinates with respect to a basis of $W$ of eigenvectors of $L(e)$, with Kovalevskaya exponents $\rho_{1},\dots,\rho_{N}$, every monomial $z^{\alpha}$ occurring in $\psi_{m}$ satisfies
\begin{equation}\label{eq:resonantmonomial}
\alpha_{1}\rho_{1}+\dots+\alpha_{N}\rho_{N}=d .
\end{equation}
In particular $I(e)=0$.
\end{lemma}

The proof below uses only the decomposition $\mathbb{A}=\mathbb{R} e\oplus e^{\perp}$ with $\theta\langle e,e\rangle =1$; diagonalizability of $L(e)$ is needed only for the reformulation \eqref{eq:resonantmonomial}.
\begin{proof}
Fix $z\in W$ and put $x=e+z$. Since $\theta\langle e,e\rangle =1$, the orthogonal projection onto $\mathbb{R} e$ is $y\mapsto\theta\langle y,e\rangle e$; let $\pi$ be the orthogonal projection onto $W$. As $ez\in W$,
\begin{equation}\label{eq:xsquare}
\begin{aligned}
x^{2}&=\bigl(1+\gamma(z)\bigr)e+2ez+\pi(z^{2}),\\
\gamma(z)&=\theta\langle z^{2},e\rangle .
\end{aligned}
\end{equation}
For $w\in W$ we have $dI_{x}(w)=d\psi_{z}(w)$, and Euler's identity $dI_{x}(x)=d\cdot I(x)$ gives
\begin{equation*}
dI_{x}(e)=d\cdot \psi(z)-d\psi_{z}(z).
\end{equation*}
Therefore, by \eqref{eq:xsquare},
\begin{equation*}
0=dI_{x}(x^{2})=\bigl(1+\gamma(z)\bigr)\bigl(d\cdot \psi(z)-d\psi_{z}(z)\bigr)+d\psi_{z}(2ez)+d\psi_{z}\bigl(\pi(z^{2})\bigr).
\end{equation*}
The terms involving $\gamma(z)$ or $\pi(z^{2})$ are of degree at least $m+1$ in $z$. Taking the homogeneous component of degree $m$ gives
\begin{equation*}
0=d\cdot \psi_{m}-d\psi_{m}(z)(z)+d\psi_{m}(z)(2ez)=d\cdot \psi_{m}-d\psi_{m}(z)(\mathcal{K}_{e}z),
\end{equation*}
which is \eqref{eq:resonance}. In Peirce coordinates, $\phi\mapsto d\phi(z)(\mathcal{K}_{e}z)$ is the derivation with $z_{i}\mapsto\rho_{i}z_{i}$, which multiplies $z^{\alpha}$ by $\sum_{i}\alpha_{i}\rho_{i}$. Comparing coefficients in \eqref{eq:resonance} gives \eqref{eq:resonantmonomial}. If $I(e)\neq0$, then $m=0$ and \eqref{eq:resonance} reads $d\cdot I(e)=0$, which is impossible since $d\geq1$.
\end{proof}

A monomial satisfying $\sum_{i}\alpha_{i}\rho_{i}=d$ is called \emph{resonant}. We apply the lemma to $\mathcal{E}$ and $\mathcal{F}$, for which $d=4$ and $d=6$, respectively.

We can now prove Theorem~\ref{mthm:C}, which we restate for convenience.

\begin{theorem}[Theorem~\ref{mthm:C}]\label{thm:rpc}
Let $\mathbb{A}$ be radial $\infty$-Hsiang. If some $e\in\operatorname{Idm}(\mathbb{A})$ satisfies $n_{-1/2}(e)=0$, then $\mathcal{E}\equiv0$. Consequently $\mathbb{A}$ is radial pseudocomposition, and $n_{-1/2}(e')=0$ for every $e'\in\operatorname{Idm}(\mathbb{A})$.
\end{theorem}
\begin{proof}
By Lemma~\ref{prop:peirce} and the assumption, the Kovalevskaya exponents on $W=e^{\perp}$ are $3$ on $\mathbb{A}_{e}(-1)$ and $0$ on $\mathbb{A}_{e}(\frac12)$. Let $\psi(z)=\mathcal{E}(e+z)$, a polynomial of degree at most $4$ on $W$. A monomial of degree $m\leq4$ with $j$ factors from $\mathbb{A}_{e}(-1)$ has
\begin{equation*}
\sum\alpha_{i}\rho_{i}-4=3j-4\in\{-4,-1,2,5,8\},
\end{equation*}
which is never zero. So no monomial is resonant, and Lemma~\ref{lem:resonance} forces $\psi\equiv0$.

Thus $\mathcal{E}$ vanishes on the affine hyperplane $e+W$. If $\langle x,e\rangle \neq0$, put $t=\theta\langle x,e\rangle $. Then $t^{-1}x-e\in W$, so $\mathcal{E}(x)=t^{4}\mathcal{E}(t^{-1}x)=0$. The set $\{x:\langle x,e\rangle \neq0\}$ is open and dense, hence $\mathcal{E}\equiv0$ by continuity. By Theorem~\ref{thm:integrals}(iii), $\mathbb{A}$ is radial pseudocomposition, and Lemma~\ref{lem:defects} gives $n_{-1/2}(e')=0$ for every $e'\in\operatorname{Idm}(\mathbb{A})$.
\end{proof}

\begin{remark}\label{rem:spectrum}
Lemma~\ref{lem:resonance} also recovers the Peirce polynomial \eqref{eq:peircepoly} using from Lemma~\ref{prop:peirce} only that $\theta\langle e,e\rangle =1$ and that $e^{\perp}$ is $L(e)$-invariant. Let $W=e^{\perp}$ and $L=L(e)|_{W}$, so that $\mathcal{K}_{e}=\operatorname{id}-2L$ on $W$. The computation of $d\mathcal{R}_{e}$ in Lemma~\ref{lem:defects} uses only $\theta\langle e,e\rangle =1$ and gives $d\mathcal{R}_{e}=M:=2L^{2}+L-\operatorname{id}$ on $W$, so that the quadratic part of $\mathcal{E}(e+z)$ is $2\langle Mz,z\rangle $. If $M\neq0$, this is the lowest part of $\mathcal{E}(e+z)$, and \eqref{eq:resonance} with $d=4$ gives $\langle Mz,(\mathcal{K}_{e}-2\operatorname{id})z\rangle =0$ for all $z\in W$. Since $M$ and $\mathcal{K}_{e}$ are commuting self-adjoint operators, $M(\mathcal{K}_{e}-2\operatorname{id})$ is self-adjoint, hence zero; that is, $(L+\operatorname{id})(2L-\operatorname{id})(2L+\operatorname{id})=0$ on $W$. If $M=0$, this identity holds trivially. Thus the Peirce eigenvalues are exactly those at which the quadratic part of the energy is either zero or resonant.
\end{remark}

\begin{remark}
Theorem~\ref{thm:rpc} explains why the absence of the eigenvalue $-\frac12$ is a global property: the energy $\mathcal{E}$ is a global object, its Hessian at any idempotent detects $n_{-1/2}$, and if this Hessian vanishes at one idempotent, then $\mathcal{E}$ vanishes everywhere. The multiplicities of the other eigenvalues are not controlled by $\mathcal{E}$, since the exponents $3$ and $0$ never combine to $4$. Example~\ref{ex:rpc} shows that they do vary.
\end{remark}

\section{Admissible cubic algebras: a resonant obstruction}\label{sec:rac}

Suppose now that $n_{-1}(e)=0$. The exponents on $W=e^{\perp}$ are then $2$ on $\mathbb{A}_{e}(-\frac12)$ and $0$ on $\mathbb{A}_{e}(\frac12)$, and the integral $\mathcal{F}$ of \eqref{eq:EF} has degree $6$. A monomial with $j$ factors from $\mathbb{A}_{e}(-\frac12)$ is resonant if and only if $2j=6$, that is $j=3$, with an arbitrary number of further factors from $\mathbb{A}_{e}(\frac12)$. Unlike for $\mathcal{E}$, resonant monomials therefore exist, and they are the only possible obstruction.

\begin{lemma}\label{prop:racobstruction}
Let $\mathbb{A}$ be radial $\infty$-Hsiang and $e\in\operatorname{Idm}(\mathbb{A})$.
\begin{itemize}
\item[{(i)}] For every $y\in\mathbb{A}_{e}(-\frac12)$, the polynomial $s\mapsto \mathcal{F}(e+sy)$ is divisible by $s^{4}$.
\item[{(ii)}] If $n_{-1}(e)=0$ and $\mathcal{F}\not\equiv0$, then the lowest-order part of $z\mapsto \mathcal{F}(e+z)$, $z\in e^{\perp}$, has degree $m\in\{4,5,6\}$ and is a combination of monomials with exactly three factors from $\mathbb{A}_{e}(-\frac12)$ and $m-3\geq1$ factors from $\mathbb{A}_{e}(\frac12)$.
\end{itemize}
\end{lemma}
\begin{proof}
(i) By Theorem~\ref{thm:hessians}, $\mathcal{F}(e)=0$, $\nabla \mathcal{F}(e)=0$, and $d^{2}\mathcal{F}_{e}$ vanishes on $\mathbb{A}_{e}(-\frac12)$, so we must show that the coefficient of $s^{3}$ vanishes. Recall that $u(w)=\langle w,w^{2}\rangle $, so that $\mathcal{F}(x)=u(x^{2})-\theta u(x)^{2}$, and note that
\begin{equation}\label{eq:Nexp}
u(e+h)=\langle e,e\rangle +3\langle e,h\rangle +3\langle eh,h\rangle +u(h).
\end{equation}
With $x=e+sy$ we have
\begin{equation*}
x^{2}=e+sh_{1}+s^{2}h_{2},\qquad h_{1}=2ey=-y,\quad h_{2}=y^{2}.
\end{equation*}
We compute the $s^{3}$-coefficient of $u(x^{2})$ from \eqref{eq:Nexp} with $h=sh_{1}+s^{2}h_{2}$. The term $3\langle e,h\rangle $ is linear in $h$ and does not contribute. The term $3\langle eh,h\rangle $ contributes
\begin{equation*}
6\langle eh_{1},h_{2}\rangle =6\langle \tfrac12y,y^{2}\rangle =3u(y),
\end{equation*}
and $u(h)$ contributes $u(h_{1})=-u(y)$. So the $s^{3}$-coefficient of $u(x^{2})$ is $2u(y)$. On the other hand,
\begin{equation*}
u(e+sy)=\langle e,e\rangle +3s\langle e,y\rangle +3s^{2}\langle ey,y\rangle +s^{3}u(y)=\langle e,e\rangle -\tfrac32s^{2}\langle y,y\rangle +s^{3}u(y),
\end{equation*}
so the $s^{3}$-coefficient of $\theta u(x)^{2}$ is $2\theta\langle e,e\rangle u(y)=2u(y)$. The two contributions cancel.

(ii) By Lemma~\ref{lem:resonance} the lowest part $\psi_{m}$ consists of resonant monomials; the values $2j-6$, $0\leq j\leq6$, vanish only for $j=3$. So $m\geq3$ and every monomial has exactly three factors from $\mathbb{A}_{e}(-\frac12)$. If $m=3$, then $\psi_{3}$ involves only coordinates on $\mathbb{A}_{e}(-\frac12)$, and $\psi_{3}(y)$ is the $s^{3}$-coefficient of $\mathcal{F}(e+sy)$, which vanishes by (i). Hence $4\leq m\leq6$.
\end{proof}

\begin{corollary}\label{cor:rac2}
Let $\mathbb{A}$ be radial $\infty$-Hsiang. If some $e\in\operatorname{Idm}(\mathbb{A})$ satisfies $n_{-1}(e)=n_{1/2}(e)=0$, then $\mathcal{F}\equiv0$, $\mathbb{A}$ is radial admissible cubic, and $n_{-1}(e')=0$ for every $e'\in\operatorname{Idm}(\mathbb{A})$.
\end{corollary}
\begin{proof}
If $\mathcal{F}\not\equiv0$, Lemma~\ref{prop:racobstruction}(ii) requires a factor from $\mathbb{A}_{e}(\frac12)=0$, which is impossible. Hence $\mathcal{F}(e+z)=0$ for $z\in e^{\perp}$, and $\mathcal{F}\equiv0$ by homogeneity, as in the proof of Theorem~\ref{thm:rpc}. Conclude with Theorem~\ref{thm:integrals}(iv) and Lemma~\ref{lem:defects}.
\end{proof}

\begin{remark}
In the remaining case, $n_{-1}(e)=0$ with $n_{-1/2}(e)\geq1$ and $n_{1/2}(e)\geq1$, Lemma~\ref{prop:racobstruction} reduces the rigidity of the eigenvalue $-1$ to the vanishing of three explicit forms, of types $(-\frac12)^{3}(\frac12)^{j}$ with $j=1,2,3$. These mixed terms are not forced to vanish by the resonance condition alone, so their vanishing when $n_{-1}(e)=0$ would have to come from the way $\mathbb{A}_{e}(\frac12)$ is acted on by the multiplication of $\mathbb{R} e\oplus\mathbb{A}_{e}(-\frac12)$, which is what the Jordan algebra structure used in the algebraic approach encodes. Since $\mathbb{A}_{e}(\frac12)$ is the tangent space of the idempotent manifold, the remaining obstruction is a statement about how the idempotent manifold bends in the directions of $\mathbb{A}_{e}(-\frac12)$.

The difference between the two rigidity statements thus has a simple numerical origin: the energy has degree $4$, which is a sum of the exponent $2$ in one way only and is not a sum of copies of $3$, whereas $\mathcal{F}$ has degree $6=3+3=2+2+2$.
\end{remark}

\section{Concluding remarks}\label{sec:final}

\medskip\noindent\textbf{Summary.}~
For a radial $\infty$-Hsiang algebra the Peirce data of an idempotent $e$ have three different interpretations. The number $n_{1/2}(e)$ is the dimension of the idempotent manifold at $e$, that is, of the maximum set of the cubic form on the sphere (Theorems~\ref{mthm:A} and \ref{mthm:MB}); dynamically, it is the number of zero Kovalevskaya exponents at the Darboux point $-e$. The numbers $n_{-1/2}(e)$ and $n_{-1}(e)$ are the ranks at $e$ of the Hessians of the two global first integrals $\mathcal{E}$ and $\mathcal{F}$ (Theorem~\ref{mthm:B}). The defining identity is equivalent to the conservation of $\mathcal{E}$, and the two classical refinements correspond to $\mathcal{E}\equiv0$ and $\mathcal{F}\equiv0$.

\medskip\noindent\textbf{Why the eigenvalues $-1$ and $-\frac12$ are global but multiplicities are not.}~
Since $\mathcal{E}$ and $\mathcal{F}$ are global, their identical vanishing is detected at any single idempotent, and for $\mathcal{E}$ this is a direct consequence of the resonance lemma (Theorem~\ref{mthm:C}). The multiplicities, on the other hand, are ranks of Hessians of the same integrals at \emph{different} points, and nothing forces them to agree. They are constant on connected components of the idempotent manifold but may change between components (Examples~\ref{ex:rpc} and \ref{ex:rac}). The same applies to the occurrence of the eigenvalue $\frac12$, which simply means that $e$ is not an isolated idempotent.

\medskip\noindent\textbf{Questions.}~
\begin{itemize}
\item[{(1)}] Is there a first integral of \eqref{eq:flow}, built from the metric and the multiplication, whose degree is not a sum of copies of $2$ and which characterizes admissible cubic algebras? Such an integral would give a resonance proof of the rigidity of the eigenvalue $-1$. Among polynomial invariants of degrees $5$ and $7$ formed from $x,x^{2},x^{3},x^{2}x^{2},xx^{3}$ and the metric, we found no first integrals other than those that vanish identically by \eqref{eq:hsiang} and \eqref{eq:V}.
\item[{(2)}] For a homogeneous potential with finitely many simple Darboux points, the Hessian eigenvalues at all Darboux points satisfy explicit global relations \cite{MaciejewskiPrzybylska05}, and similar relations hold for the multipliers at the fixed points of a polynomial endomorphism of projective space \cite{Guillot04}, in particular for the spectra of $L(c)$ over the idempotents $c$ of a generic commutative algebra \cite{KrasnovTkachev19}. For radial $\infty$-Hsiang algebras the Darboux points are in general not isolated, since $n_{1/2}(e)$ is the dimension of the family. Is there a version of these relations, in the spirit of residue formulas for nonisolated singular sets, which relates the Peirce numbers of the different components of $\operatorname{Idm}(\mathbb{A})$? In both Examples~\ref{ex:rpc} and \ref{ex:rac}, the isolated idempotent and the positive-dimensional component have quite different Peirce numbers.
\item[{(3)}] Is every connected component of $\operatorname{Idm}(\mathbb{A})$ an orbit of the automorphism group, as in all known examples? In Example~\ref{ex:rpc} the sphere of idempotents is an orbit of the orthogonal group $O(m)$ acting by $(r,x)\mapsto(r,gx)$, and in Example~\ref{ex:rac} the projective plane is an orbit of $SO(3)$ acting by $x\mapsto gxg^{T}$.
\item[{(4)}] The local part of Theorem~\ref{mthm:A} does not use positivity of the metric (Remark~\ref{rem:indefinite}). For an algebra with a nondegenerate invariant metric over a general field of characteristic not dividing $30$, is the scheme of idempotents smooth of relative dimension $n_{1/2}$ at every idempotent?
\end{itemize}

\section*{Acknowledgments}

\noindent
The second author is grateful to Yakov Krasnov for many discussions over the years, and in particular for the idea of using the dictionary between Darboux points and Kovalevskaya exponents of a quadratic system, on the one hand, and idempotents and their Peirce spectra, on the other. Some of the ideas of this paper go back in a preliminary form to these discussions around \cite{KrasnovTkachev19} and to our joint unfinished projects on quadratic differential equations.
The second author also thanks Evgeniy Lokharu for stimulating discussions on the Lyapunov--Schmidt reduction in bifurcation theory, and in particular on the question of when the bifurcation equations are satisfied identically; these discussions helped to shape the argument in Section~\ref{sec:euclid}.


\end{document}